\pdfoutput=1

\documentclass[11pt]{amsart} 

\usepackage[margin=2.9cm]{geometry}	
\usepackage{amsthm}	
\usepackage{amsmath}	
\usepackage{amssymb}	
\usepackage{extarrows}	
\usepackage{mathtools}	
\usepackage{enumerate}	
\usepackage{tabularx}	
\usepackage{xcolor}	
\usepackage{tikz}	
\usepackage{relsize}	
\usepackage{multicol}	
\usepackage{graphicx}	
\usepackage[font={small}]{caption}	
\usepackage{subcaption}	
\usepackage{wrapfig}	
\usepackage{marginnote}	
\usepackage{comment}	
\usepackage[normalem]{ulem}	
\usepackage{subfiles}	
\usepackage{xifthen}	
\usepackage{xparse}	
\usepackage{etoolbox}	
\usepackage{environ}	

\usepackage[LGR,T1]{fontenc}	
\usepackage[utf8]{inputenc}	
\usepackage{lmodern}	
\usepackage[greek,english]{babel}	
\usepackage{alphabeta}	
\usepackage{csquotes}	

\usepackage[doi=true,isbn=true,giveninits=true,maxbibnames=99,backend=biber]{biblatex}	
\usepackage{xurl}	
\usepackage{cleveref}	
\hypersetup{colorlinks=true,linktoc=section,urlcolor={blue!80!black},citecolor={green!60!black}}
\nocite{
		Tao2008,
		Tol2006,
}	
\makeatletter	
\AtBeginDocument{
\patchcmd\label@noarg{\edef\@tempb}{\protected@edef\@tempb}{}{}}	
\makeatother	

\newtheorem{thm}{Theorem}
\newtheorem{lem}[thm]{Lemma}	
\newtheorem*{thm*}{Theorem}	
\newtheorem*{lem*}{Lemma}	
\newtheorem*{prop*}{Proposition}	
\newtheorem*{cor*}{Corollary}	
\newtheorem*{dfn*}{Definition}	
\newtheorem*{rem*}{Remark}	
\newtheorem*{rems*}{Remarks}	
\newtheorem*{que*}{Question}	
\newtheorem*{not*}{Notations}	
\newtheorem*{eg*}{Example}	
\numberwithin{equation}{section}	

		\def \P{{\mathbb{P}}}		
\def \C{{\mathbb{C}}}				
\def \D{{\mathbb{D}}}		\def \R{{\mathbb{R}}}		
\def \E{{\mathbb{E}}}				
	\def \N{{\mathbb{N}}}								
\def \cD{{\mathcal{D}}}	\def \cK{{\mathcal{K}}}	\def \cR{{\mathcal{R}}}		
\def \cE{{\mathcal{E}}}				
\def \cF{{\mathcal{F}}}		\def \cT{{\mathcal{T}}}							
\def \cG{{\mathcal{G}}}									
\NewCommandCopy{\oldin}{\in}	
\DeclareRobustCommand{\in}{\oldin\nolinebreak[4]}												
\makeatletter	
\def\hrulefill{\noindent\leavevmode\leaders\hrule\hfill\kern\z@}								
\makeatother	
\def \such{\,:\,}	
\newcommand{\inline}[1]{\quad\text{#1}\quad}	
\newcommand{\afterline}[1]{\qquad\text{#1~}}	
\DeclarePairedDelimiterX{\inner}[2]{\langle}{\rangle}{#1,#2}	

\def\?{\textcolor{red}{\textbf{check/why?}}}	

\DeclareMathOperator{\proj}{proj}	
\DeclareMathOperator{\Fav}{Fav}	
\newcommand{\needle}{l^\perp}	
\newcommand{\flipup}[1]{\ifnumodd{#1}{node}{node[yshift=12pt]}}
\usetikzlibrary{
	arrows.meta,
	math,
	decorations.pathreplacing,
	decorations.pathmorphing,
	calligraphy}	
\newlength{\dotlength}	
\makeatletter	
\tikzset{
    dot diameter/.store in=\dot@diameter,
    dot diameter=0.5\dotlength,
    dot spacing/.store in=\dot@spacing,
    dot spacing=0.8em,
    dots/.style={
        line width=\dot@diameter,
        line cap=round,
        dash pattern=on 0pt off \dot@spacing
    }
}
\makeatother	

\title{The Buffon's needle problem for random planar disk-like Cantor sets}
\author[D.~Vardakis]{Dimitris Vardakis}
\address{Michigan State University}
\email{jimvardakis@gmail.com}
\author[A.~Volberg]{Alexander Volberg}
\address{Michigan State University}
\email{volberg@msu.edu}
\thanks{AV is supported by NSF grants DMS 1900286 and DMS 2154402}
\keywords{Favard length, Buffon needle, Hausdorff measure, random Cantor sets}

\begin{document}

\begin{abstract} 
	We consider a model of randomness for self-similar Cantor sets of finite and positive $1$-Hausdorff measure.
	We find the sharp rate of decay of the probability that a Buffon needle lands $\delta$-close to a Cantor set of this particular randomness.
	Two quite different models of randomness for Cantor sets, by Peres and Solomyak, and by Shiwen Zhang, appear to have the same order of decay for the Buffon needle probability: $\frac{c}{\log\frac{1}{\delta}}$.
	In this note, we prove the same rate of decay for a third model of randomness, which asserts a vague feeling that any ``reasonable'' random Cantor set of positive and finite length will have Favard length of order $\frac{c}{\log\frac{1}{\delta}}$ for its $\delta$-neighbourhood.
	The estimate from below was obtained long ago by Mattila.
\end{abstract}

\maketitle 
\tableofcontents 

\section{Introduction} \label{sec:intro}

Let $E$ be a subset of the unit disk, $\D$. The \emph{Buffon needle problem} wants to determine the probability with which a random needle or line intersects $E$ provided that it already intersects the unit disk. At the same time, let $l_\theta$ be the line passing through the origin and forming angle $\theta$ with the horizontal axis. The \emph{Favard length} of $E$ is the average length of the projection of $E$ onto $l_\theta$ when averaging over all angles $\theta$. It turns out these two quantities are proportional.

\medskip

Now, consider the following picture: let us have $L$ many ($L\geq3$) disjoint closed disks $(D_1,\dots,D_L)$ of diameter $1/L$ and strictly inside $\D$. These are disks of the first generation. Consider also a piecewise affine map $f=(f_1,\dots,f_L)$ from those disks onto $\D$. Then, $f^{-1}(\D)=D_1\cup\dots\cup D_L$. Furthermore, $f^{-1}(D_1\cup\dots\cup D_L)$ it consists of $L^2$ disks (groups of $L$ many disks in each $D_i$); we call those disks of the second generation. We can iterate this procedure: denoting by $U_n$ the union of disks of the $n$-th generation, where $U_1:=D_1\cup\dots\cup D_L$, we form the self-similar Cantor set $\cK=\bigcap_{n=1}^\infty U_n$. This has positive and finite $1$-dimensional Hausdorff measure; thus it is completely unrectifiable in the sense of Besicovitch \cite{Mat2015}; and thus its Favard length is zero \cite{Mat2015}.

Of course, the disks can be replaced by other shapes. For example, $U_1$ can consist of $L$ disjoint squares with side-length $1/L$ inside the unit square $[0,1]^2$ (where the word \enquote{strictly} can be omitted but \enquote{disjoint} cannot). One of such Cantor sets is a rather \enquote{famous}, namely the $1/4$-corner Cantor set, $\cK_{1/4}$ (see \cite{Mat2004}).

\medskip

The $L^{-n}$-neighbourhood of such sets is roughly $U_n$, and therefore its Favard length
\[\Fav(U_n)\to 0,\quad\text{as }n\to\infty.\]
But what is $\Fav(U_n)$, or what is the speed with which $\Fav(U_n)$ decreases? Nobody knows exactly, but there has been considerable interest in recent years. It is now clear that the answer may depend on several factors; the magnitude of $L$; the geometry of $U_1$; the subtle algebraic and number theoretic properties of a certain trigonometric sum built by the centres of the disks of the first generation. See \cite{BonLabVol2014,BonVol2010,BonVol2012,LabZha2010,NazPerVol2011} and the survey paper \cite{Lab2014}.

For the $1/4$-corner Cantor set $\cK_{1/4}$ in particular, the best known estimate from above for its $4^{-n}$-neighbourhood is
\begin{equation*} \label{eq:16}
	\Fav(N_{4^{-n}}(\cK_{1/4}))\leq\frac{C_\varepsilon}{n^{\frac{1}{6}-\varepsilon}},\quad\forall \varepsilon>0,
\end{equation*}
for all large $n$. We suspect that this estimate can be improved to
\begin{equation*} \label{eq:1eps}
	\Fav(N_{4^{-n}}(\cK_{1/4}))\leq\frac{C_\varepsilon}{n^{1+\varepsilon}},\quad\forall \varepsilon>0,
\end{equation*}
but at this moment this is only a conjecture.

\medskip

On the other hand, there is a universal estimate from below obtained in \cite{Mat1975} for every self-similar Cantor set constructed as above:
\begin{equation} \label{below}
	\Fav (N_{4^{-n}}(\cK))\geq\frac{c}{n}.
\end{equation}
For any concrete set, this bound from below could be improved. In fact, it is proven in \cite{BatVol2010} that for the same $1/4$-corner Cantor set $\cK_{1/4}$
\[\Fav(N_{4^{-n}}(\cK_{1/4}))\geq\frac{c\log n}{n}.\]

\medskip

For random Cantor sets the situation should be simpler. With large probability, Mattila's lower estimate \eqref{below} is met by \emph{the same} estimate from above (with a different constant). The problem is that in general there can be many different models of randomness.

\medskip

In this note, we are interested in an analogue of the random Cantor set appearing in \cite{PerSol2002} and in \cite{Zha2019}. In our case, this will come from the random Cantor disks constructed below at \Cref{sec:the_disks}. The model of randomness presented here is somewhat different from the ones in the above two papers, but it amazingly exhibits the same behaviour, as we'll see below in our main \Cref{thm:Favard}, which we contrast with \cite[Theorem 2.2]{PerSol2002} and \cite[Theorem 1]{Zha2019}.

In particular, we prove an analogue of \cite[Theorem 1]{Zha2019}. Unfortunately, the randomness of the disk model we study here is not equivalent to that of the random (square) Cantor set $\cR=\bigcap_{n=0}^\infty\cR_n$ from \cite{PerSol2002}, but it is nonetheless closer compared to the one constructed in \cite{Zha2019}. The essential difference between \cite{Zha2019} and our consideration are the angles $\omega_n^1,\omega_n^2,\dots,\omega_n^{4^{n-1}}$, which are here allowed to be distinct and independent whereas in \cite{Zha2019} are all equal. So, our model is a little \enquote{more random} than the random Cantor sets of Zhang in \cite{Zha2019}.

We introduce our notations ---some borrowed from \cite{Zha2019}--- in the next \Cref{sec:the_disks}. The problem of interest, namely the \emph{Favard length} of a random planar disk-like Cantor set, is explained in \Cref{sec:Favard_length}. Our results and their proofs are postponed to \Cref{sec:the_main_lemma,sec:proof_of_lemma}. In \Cref{sec:compa}, we compare the differences and difficulties between our work and that of Peres and Solomyak's and Zhang's.

\section{Cantor Disks} \label{sec:the_disks}

Our work will be heavy on notation; without any ado let us introduce our basic \enquote{vocabulary}.

\medskip

The letter $n$ will stand for a (large) positive integer.

The letter $\omega$ will be used to denote angles with values inside the interval $[0,\frac{\pi}{2}]$. Now, let us consider a word of length $n$ made of the alphabet of angles in $[0,\frac{\pi}{2}]$, i.e. a word of the form $\omega_1\omega_2\cdots \omega_n$. The subscript in $\omega_k$ denotes the position of the angle $\omega_k$ within such a word of length $n$. We refer to the position of an angle within a word as the \emph{depth} of that angle.

Our operators, which we will introduce below, are such that every choice of an angle, say, $\omega_1$ necessitates four different independent choices for the angle $\omega_2$; every choice of the angle $\omega_2$ necessitates four different independent choices for the angle $\omega_3$; and so on up until depth $n$ where we will have $4^{n-1}$ different angles $\omega_n$. In order to differentiate between all those, for each $j_k=1,2,\dots,4^{k-1}$ we write $\omega_k^{j_k}$ for the $j_k$-th choice of an angle $\omega_k$ at depth $k$. Notice there are $4^{k-1}$ such choices. Therefore, a typical word from our alphabet of angles looks as follows, where we note that $\omega_k^{j_k}\in[0,\frac{\pi}{2}]$:
\begin{align*}	\label{def:words_of_angles}
	\omega_1^{j_1}\omega_2^{j_2}\!\cdots \omega_k^{j_k}\!\cdots \omega_n^{j_n}\afterline{where}
	\begin{array}{r@{~}lr@{~}l}
		j_1	&	=1\,,		&	j_k	&	=1,2,\dots,4^{k-1}\,,\\
		j_2	&	=1,2,3,4\,,	&		&	\cdots\\
		&	\cdots		&	j_n	&	=1,2,\dots,4^{n-1}.
	\end{array}
\end{align*}

At certain instances, we need to consider cumulatively all angles of a certain depth; given a collection of words of length $n$, for each $k=1,2,\dots,n$ let $\omega_k'$ be the collection of all $4^{k-1}$ many angles at depth $k$, that is
\begin{equation*}	\label{def:angles_of_depth}
	\omega_k'=(\omega_k^1,\dots,\omega_k^{4^{k-1}}).
\end{equation*}
With this notation, we may use the symbols $\omega_1$, $\omega_1^{j_1}$, $\omega_1^1$, and $\omega_1'$ interchangeably as these all refer to the same single angle.

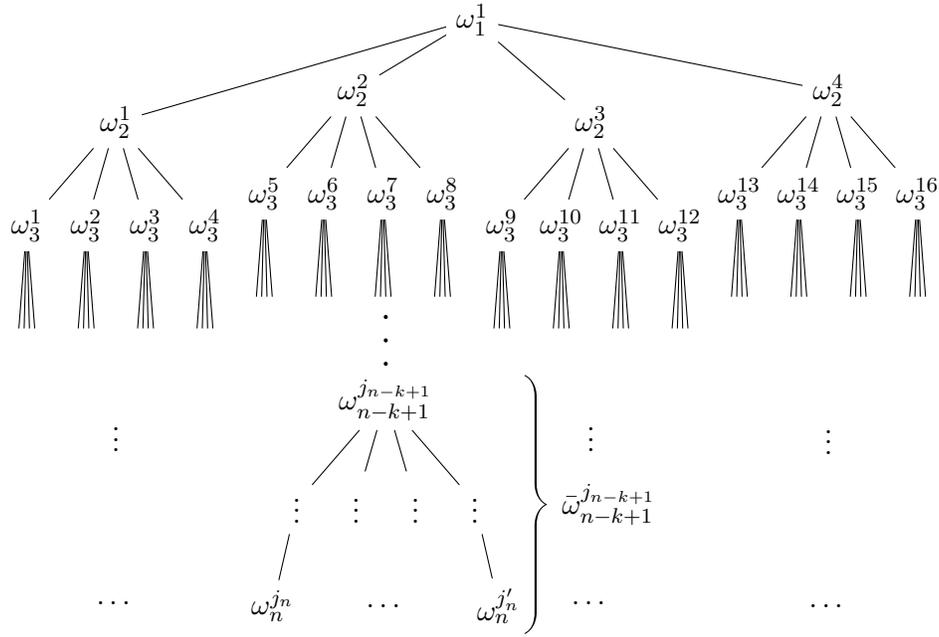
\begin{figure}
	\centering
	\begin{tikzpicture}[%
		level distance={3\baselineskip},
		level 1/.style={sibling distance=0.20\textwidth},
		level 2/.style={sibling distance=0.05\textwidth},
		level 3/.style={sibling distance=2pt},
		]
		\node (start) {$\omega_1^1$}
		child foreach \v in {1,...,4} {\flipup{\v} {$\omega_2^\v$}
			child foreach \x in {1,...,4} {node {\pgfmathparse{int(4*(\v-1)+\x)} $\omega_3^{\pgfmathresult}$}
				child foreach \y in {0,...,3}
			}
		};
		\node at (start-2-3-3) {\phantom{M}}
		[level 3/.style={sibling distance=0.04\textwidth}]
			child {node (R) {$\omega_{n-k+1}^{j_{n-k+1}}$} edge from parent[dots]
				child {node {$\vdots$} edge from parent[solid,thin]
					child {node (L1) {$\omega_n^{j_n}$}
					}
					child[missing]
				}
				child foreach \z in {2,3} {node {$\vdots$} edge from parent[solid,thin]}
				child {node {$\vdots$} edge from parent[solid,thin]
					child[missing]
					child {node (L2) {$\omega_n^{j_n'}$}
						child[grow=north,yshift=4\baselineskip] {node (R2) {} edge from parent[draw=none]}
					}
				}
			};
		\path (start-1-2-4) -- (start-1-3-1) node[midway,yshift=-3\baselineskip] {$\vdots$}
			child {node[yshift=-2\baselineskip] {$\dots$} edge from parent[draw=none]};
		\path (start-3-2-4) -- (start-3-3-1) node[midway,yshift=-3\baselineskip] {$\vdots$}
			child {node[yshift=-2\baselineskip] {$\dots$} edge from parent[draw=none]};
		\path (start-4-2-4) -- (start-4-3-1) node[midway,yshift=-4\baselineskip] {$\vdots$}
			child {node[yshift=-2\baselineskip] {$\dots$} edge from parent[draw=none]};
		\path (L1) -- (L2) node [midway] {$\dots$};
		\draw [thick,decorate,decoration={calligraphic brace,raise=10pt,amplitude=8pt}] (R2) -- (L2.south) node[midway,right=20pt] {$\bar{\omega}_{n-k+1}^{j_{n-k+1}}$};
	\end{tikzpicture}
	\caption{The $n$-tree of angles with root $\omega_1^1$, all angles of depth 1 to 3, and the subtree $\bar{\omega}_{n-k+1}^{j_{n-k+1}}$ of height $k$ with root $\omega_{n-k+1}^{j_{n-k+1}}$ at depth $n-k+1$. All angles are enumerated. (The exact values of $j_{n}$ and $j_{n}'$ depend on $j_{n-k+1}$.)}
	\label{fig:trees}
\end{figure}

All the above give to our angles the structure of a rooted tree of height $n$ with root $\omega_1$ and such that each parent has four children as in \Cref{fig:trees}. The vertexes have values in $[0,\frac{\pi}{2}]$, and are independent from each other and from their predecessors and ancestors. This tree we denote by $\omega_1'\cdots \omega_n'$; the trimmed tree with root $\omega_1$ and height $k$ we denote as $\omega_1'\cdots \omega_k'$ (for any $k=1,2,\dots,n$). For the subtree of height $n-k+1$ with root $\omega_k^{j_k}$, which reaches up to the leaves (that is, from depth $k$ till depth $n$ with starting vertex $\omega_k^{j_k}$) we write $\bar{\omega}_k^{j_k}$. Later on, we will be working with rooted subtrees of the form $\bar{\omega}_{n-k+1}^{j_{n-k+1}}$. To reiterate, $\bar{\omega}_{n-k+1}^{j_{n-k+1}}$ consists of the angle $\omega_{n-k+1}^{j_{n-k+1}}$ (as its root) along with all the angles from depth $n-k+1$ till depth $n$ (which have $\omega_{n-k+1}^{j_{n-k+1}}$ as an ancestor). This has height $k$. Alternatively, $\bar{\omega}_{n-k+1}^{j_{n-k+1}}$ is the collection of all the words (from our alphabet of angles) which have depth (i.e. length) $k$ and the first letter is $\omega_{n-k+1}^{j_{n-k+1}}$. There are $4^{n-k}$ such words.

\bigskip

Next, we will need to introduce certain operators and sets. The main objects of interest will be the operators $\cD_k$ ($k=0,1,\dots,n$) which will act on trees of angles of depth $k$. To understand these we need some auxiliary constructions first.

\medskip

For any angle $\omega$ and for $\alpha=0,1,2,3$ consider the transformations
\begin{equation}	\label{def:shrinking}
	T_\alpha^\omega(z)=\frac{1}{4}z+\frac{3}{4}e^{(\alpha\frac{\pi}{2}-\omega)i}
\end{equation}
where $z$ is any number on the complex plane $\C$. Observe that if $\D$ is the unit disk, $T_0^0(\D)$, $T_1^0(\D)$, $T_2^0(\D)$, and $T_3^0(\D)$ are disks of radius $1/4$ centred respectively at $(3/4,0)$, $(0,3/4)$, $(-3/4,0)$, and $(0,-3/4)$. Introducing an angle $\omega$ in $T_\alpha^\omega(\D)$, rotates (about $(0,0)$) the aforementioned disks by angle $\omega$ in the clockwise direction.

Moreover, given an angle $\omega_k^{j_k}$ from depth $k$ let $\Omega_k^{j_k}$ be the set
\begin{equation*}	\label{def:aux_rotated_sets}
	\Omega_k^{j_k}=\bigcup_{\alpha=0}^3\frac{1}{4^{k-1}}T_\alpha^{\omega_k^{j_k}}(\D).
\end{equation*}
That is, $\Omega_k^{j_k}$ is a collection of four disks of radius $4^{-k}$ with centres $(0,\pm3/4^k)$ and $(\pm3/4^k,0)$ rotated clockwise by $\omega_k^{j_k}$. An example of $\Omega_2^1$ appears on \Cref{fig:disks}.

We also give an enumeration to all the disks for all depths. We number the disks of $\Omega_k^{j_k}$ so that $\frac{1}{4^{k-1}}T_\alpha^{\omega_k^{j_k}}(\D)$ is the $(4j_k-3+\alpha)$-th disk at depth $k$. We call this the \emph{$k$-depth enumeration} (of the disks lying at depth $k$). Illustratively, we note $\frac{1}{4^{k-1}}T_0^{\omega_k^1}(\D)$, $\frac{1}{4^{k-1}}T_1^{\omega_k^1}(\D)$, $\frac{1}{4^{k-1}}T_0^{\omega_k^{4^{k-1}}}(\D)$, $\frac{1}{4^{k-1}}T_3^{\omega_k^{4^{k-1}}}(\D)$ are the 1st, 2nd, $(4^k-3)$-th, $4^k$-th disks of depth $k$. We retain this enumeration as we translate these disks at different positions on the plane. This will be useful to track down each disk at each step so that our subsequent constructions make better sense.

\subsection{The \texorpdfstring{$\cD$}{D} operators}	\label{def:the_oprators}

Now, we are ready to introduce our main protagonists. The operator $\cD_k$ acts on the collection of trees (of angles) of height $k$ and for each such tree outputs a certain collection of $4^k$ disks of radius $4^{-k}$. We define these inductively below.

\medskip

To begin with, set $\cD_0=\D$ to be the unit disk.

Next, we define $\cD_1$ by
\begin{equation}	\label{eq:first_iteration}
	\cD_1(\omega_1')=\Omega_1^1=\bigcup_{\alpha=0}^3T_\alpha^{\omega_1^1}(\cD_0),
\end{equation}
that is, $\cD_1(\omega_1')$ consists of four disks of radius $1/4$ centred at $(0,\pm3/4)$ and $(\pm3/4,0)$ rotated clockwise by $\omega_1$. Recall these disks are enumerated as in $\Omega_1^1$ above.

For the operator $\cD_2$, consider a tree of height $2$, $\omega_1'\omega_2'$, which consists of the angles $\omega_1^1$, and $\omega_2^1,\omega_2^2,\omega_2^3,\omega_2^4$. Then, we define $\cD_2(\omega_1'\omega_2')$ to be the collection of disks constructed as follows: Replace the 1st, 2nd, 3rd and 4th disk of $\cD_1(\omega_1')$ respectively by $\Omega_2^1$, $\Omega_2^2$, $\Omega_2^3$ and $\Omega_2^4$. By ``replacing'' we mean the translation of $\Omega_2^j$ in such a way that $(0,0)$ is translated to the centre of the $j$-th disk of $\cD_1(\omega_1')$.

\begin{figure}[h]
	\centering
	\begin{tikzpicture}[scale=4.5]
		\color{cyan}
		\draw (0,0) circle (1);
		\node at (-4/6,6/7) {$\cD_0$};
		
		\color{blue}
		\node at (-5/9,-3/9) {$\cD_1(\omega_1')$};
		\draw [opacity=0.1] (-1/2,0) -- (1/2,0);
		\draw [opacity=0.1] (0,-1/2) -- (0,1/2);
		\draw [opacity=0.1] (0, 3/4) circle (1/4);
		\draw [opacity=0.1] (-3/4,0) circle (1/4);
		\draw [opacity=0.1] (3/4, 0) circle (1/4);
		\draw [opacity=0.1] (0,-3/4) circle (1/4);
		\draw [opacity=0.4] ({1/2*cos(60)},	{1/2*sin(60)}) -- ({1/2*cos(240)},{1/2*sin(240)});
		\draw [opacity=0.4] ({1/2*cos(150)},{1/2*sin(150)}) -- ({1/2*cos(-30)},{1/2*sin(-30)});
		\draw ({3/4*cos(60)},	{3/4*sin(60)}) circle (1/4);
		\draw ({3/4*cos(150)},	{3/4*sin(150)}) circle (1/4);
		\draw ({3/4*cos(240)},	{3/4*sin(240)}) circle (1/4);
		\draw ({3/4*cos(330)},	{3/4*sin(330)}) circle (1/4);
		\node at (1/10,4/10) {$\omega_1$};
		\draw [opacity=0.4] (0,0) -- (0,1/3) arc (90:60:1/3);
		
		\color{black}
		\node at (6/9,7/20) {$\cD_2(\omega_1'\omega_2')$};
		%
		\node at ({-1/25+3/4*cos(330)},{1/12+3/4*sin(330)}) {$\scriptstyle \omega_2^1$};
		\draw [opacity=0.4] ({3/4*cos(330)},{3/4*sin(330)}) -- ++(90:{1/15}) arc (90:55:{1/15});
		\draw [opacity=0.4] ({1/8*cos(55)}, {1/8*sin(55)}) +({3/4*cos(330)},{3/4*sin(330)}) -- ({1/8*cos(235)+3/4*cos(330)},{1/8*sin(235)+3/4*sin(330)});
		\draw [opacity=0.4] ({1/8*cos(145)},{1/8*sin(145)})+({3/4*cos(330)},{3/4*sin(330)}) -- ({1/8*cos(325)+3/4*cos(330)},{1/8*sin(325)+3/4*sin(330)});
		\draw ({3/16*cos(55)}, {3/16*sin(55)}) +({3/4*cos(330)},{3/4*sin(330)}) circle (1/4^2);
		\draw ({3/16*cos(145)},{3/16*sin(145)})+({3/4*cos(330)},{3/4*sin(330)}) circle (1/4^2);
		\draw ({3/16*cos(235)},{3/16*sin(235)})+({3/4*cos(330)},{3/4*sin(330)}) circle (1/4^2);
		\draw ({3/16*cos(325)},{3/16*sin(325)})+({3/4*cos(330)},{3/4*sin(330)}) circle (1/4^2);
		%
		\node at ({1/18+3/4*cos(60)}, {1/10+3/4*sin(60)}) {$\scriptstyle \omega_2^2$};
		\draw [opacity=0.4] ({3/4*cos(60)},{3/4*sin(60)}) -- ++(90:1/15) arc (90:30:1/15);
		\draw [opacity=0.1] ({-1/8+3/4*cos(60)},{3/4*sin(60)}) -- ({1/8+3/4*cos(60)},{3/4*sin(60)});
		\draw [opacity=0.1] ({3/4*cos(60)}, {3/4*sin(60)-1/8}) -- ({3/4*cos(60)},{3/4*sin(60)+1/8});
		\draw [opacity=0.1] ({3/4*cos(60)},{3/4*sin(60)})+(1/4-1/16, 0) circle (1/4^2);
		\draw [opacity=0.1] ({3/4*cos(60)},{3/4*sin(60)})+(-1/4+1/16,0) circle (1/4^2);
		\draw [opacity=0.1] ({3/4*cos(60)},{3/4*sin(60)})+(0, 1/4-1/16) circle (1/4^2);
		\draw [opacity=0.1] ({3/4*cos(60)},{3/4*sin(60)})+(0,-1/4+1/16) circle (1/4^2);
		\draw [opacity=0.4] ({1/8*cos(30)},{1/8*sin(30)}) +({3/4*cos(60)}, {3/4*sin(60)}) -- ({1/8*cos(210)+3/4*cos(60)},{1/8*sin(210)+3/4*sin(60)});
		\draw [opacity=0.4] ({1/8*cos(120)},{1/8*sin(120)})+({3/4*cos(60)},{3/4*sin(60)}) -- ({1/8*cos(-60)+3/4*cos(60)},{1/8*sin(-60)+3/4*sin(60)});
		\draw ({3/16*cos(30)}, {3/16*sin(30)}) +({3/4*cos(60)},{3/4*sin(60)}) circle (1/4^2);
		\draw ({3/16*cos(120)},{3/16*sin(120)})+({3/4*cos(60)},{3/4*sin(60)}) circle (1/4^2);
		\draw ({3/16*cos(210)},{3/16*sin(210)})+({3/4*cos(60)},{3/4*sin(60)}) circle (1/4^2);
		\draw ({3/16*cos(300)},{3/16*sin(300)})+({3/4*cos(60)},{3/4*sin(60)}) circle (1/4^2);
		%
		\node at ({1/18+3/4*cos(150)},{1/10+3/4*sin(150)}) {$\scriptstyle \omega_2^3$};
		\draw [opacity=0.4] ({3/4*cos(150)},{3/4*sin(150)}) -- ++(90:{1/15}) arc (90:10:{1/15});
		\draw [opacity=0.4] ({1/8*cos(10)}, {1/8*sin(10)}) +({3/4*cos(150)},{3/4*sin(150)}) -- ({1/8*cos(190)+3/4*cos(150)},{1/8*sin(190)+3/4*sin(150)});
		\draw [opacity=0.4] ({1/8*cos(100)},{1/8*sin(100)})+({3/4*cos(150)},{3/4*sin(150)}) -- ({1/8*cos(-80)+3/4*cos(150)},{1/8*sin(-80)+3/4*sin(150)});
		\draw ({3/16*cos(10)}, {3/16*sin(10)}) +({3/4*cos(150)},{3/4*sin(150)}) circle (1/4^2);
		\draw ({3/16*cos(100)},{3/16*sin(100)})+({3/4*cos(150)},{3/4*sin(150)}) circle (1/4^2);
		\draw ({3/16*cos(190)},{3/16*sin(190)})+({3/4*cos(150)},{3/4*sin(150)}) circle (1/4^2);
		\draw ({3/16*cos(280)},{3/16*sin(280)})+({3/4*cos(150)},{3/4*sin(150)}) circle (1/4^2);
		%
		\node at ({1/12+3/4*cos(240)},{1/30+3/4*sin(240)}) {$\scriptstyle \omega_2^4$};
		\draw [opacity=0.4] ({3/4*cos(240)},{3/4*sin(240)}) -- ++(90:{1/15}) arc (90:65:{1/15});
		\draw [opacity=0.4] ({1/8*cos(65)}, {1/8*sin(65)}) +({3/4*cos(240)},{3/4*sin(240)}) -- ({1/8*cos(245)+3/4*cos(240)},{1/8*sin(245)+3/4*sin(240)});
		\draw [opacity=0.4] ({1/8*cos(155)},{1/8*sin(155)})+({3/4*cos(240)},{3/4*sin(240)}) -- ({1/8*cos(335)+3/4*cos(240)},{1/8*sin(335)+3/4*sin(240)});
		\draw ({3/16*cos(65)}, {3/16*sin(65)}) +({3/4*cos(240)},{3/4*sin(240)}) circle (1/4^2);
		\draw ({3/16*cos(155)},{3/16*sin(155)})+({3/4*cos(240)},{3/4*sin(240)}) circle (1/4^2);
		\draw ({3/16*cos(245)},{3/16*sin(245)})+({3/4*cos(240)},{3/4*sin(240)}) circle (1/4^2);
		\draw ({3/16*cos(335)},{3/16*sin(335)})+({3/4*cos(240)},{3/4*sin(240)}) circle (1/4^2);
		
		\color{teal}
		\node at (4/18,1/16) {$\Omega_2^1$};
		\draw [opacity=0.2] (0,0) circle (1/4);
		\draw [-stealth] ({9/32*cos(-20)}, {9/32*sin(-20)}) -- ({9/32*cos(-20)+6/32*cos(-30)}, {9/32*sin(-20)+6/32*sin(-30)});
		\node at (-1/25,1/10) {$\scriptstyle \omega_2^1$};
		\draw [opacity=0.2] (0,0) -- (0,1/10) arc (90:55:1/10);
		\draw [opacity=0.4] ({1/8*cos(55)}, {1/8*sin(55)})	-- ({1/8*cos(235)},{1/8*sin(235)});
		\draw [opacity=0.4] ({1/8*cos(145)},{1/8*sin(145)})+({3/4*cos(330)},{3/4*sin(330)}) -- ({1/8*cos(325)+3/4*cos(330)},{1/8*sin(325)+3/4*sin(330)});
		\draw ({3/16*cos(55)}, {3/16*sin(55)})	circle (1/4^2);
		\draw ({3/16*cos(145)},{3/16*sin(145)})	circle (1/4^2);
		\draw ({3/16*cos(235)},{3/16*sin(235)})	circle (1/4^2);
		\draw ({3/16*cos(325)},{3/16*sin(325)})	circle (1/4^2);
	\end{tikzpicture}
	\caption{The collections \textcolor{cyan}{$\cD_0$}, \textcolor{blue}{$\cD_1(\omega_1')$}, $\cD_2(\omega_1'\omega_2')$ and \textcolor{teal}{$\Omega_2^1$}.}
	\label{fig:disks}
\end{figure}
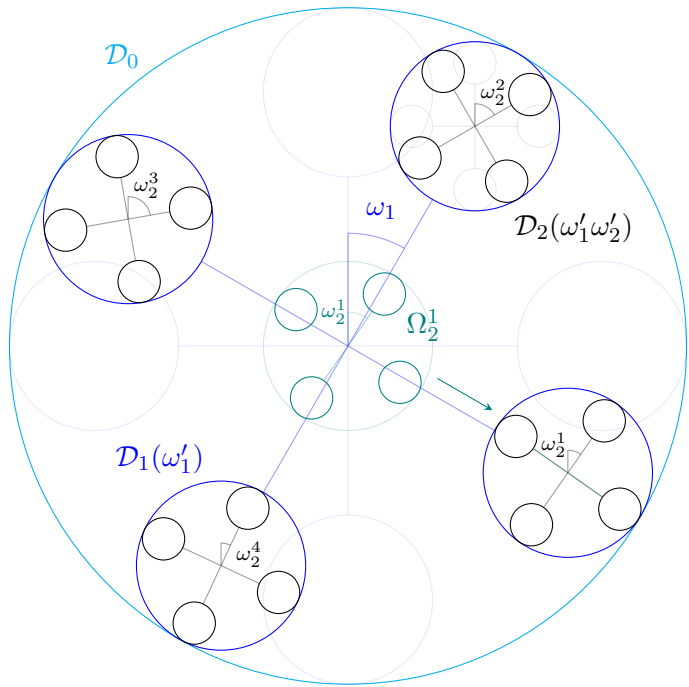

Consequently, $\cD_2(\omega_1'\omega_2')$ consists of $4^2$ disks of radius $4^{-2}$ translated appropriately so that each $\Omega_2^{j_2}$ (which is a collection of four disks) replaces one of the disks from $\cD_1(\omega_1')$. For example, the set $\Omega_2^1$ is in fact a subset of the 1st disk of $\cD_1(\omega')$; actually $\cD_2(\omega_1'\omega_2')\subset\cD_1(\omega_1')$. Again, the disks comprising $\cD_2(\omega_1'\omega_2')$ are enumerated to match $\Omega_2^1$, $\Omega_2^2$, $\Omega_2^3$ and $\Omega_2^4$ as we described above. Also see \Cref{fig:disks}.

Continuing inductively, the operator $\cD_k$ acts on the tree $\omega_1'\cdots \omega_k'$ in this manner: Consider the collection $\cD_{k-1}(\omega_1'\cdots \omega_{k-1}')$. These are $4^{k-1}$ many (enumerated) disks. Replace the 1st of them by $\Omega_k^1$, the 2nd of them by $\Omega_k^2$, etc., until every disk of $\cD_{k-1}(\omega_1'\cdots \omega_{k-1}')$ has been replaced by four smaller ones. This replacement is done so that $(0,0)$, as the \enquote{centre} of $\Omega_k^j$, is translated to the centre of the $j$-th disk of $\cD_{k-1}(\omega_1'\cdots \omega_{k-1}')$. That is, we substitute the $j$-th disk (from depth $k-1$) with the $(4j-3)$-, $(4j-2)$-, $(4j-1)$-, and $4j$-th disks of depth $k$. The resulting collection, which has $4^k$ many disks of radius $4^{-k}$, is $\cD_k(\omega_1'\dots \omega_k')$. It holds that $\cD_k(\omega_1'\cdots \omega_k')\subset\cD_{k-1}(\omega_1'\cdots \omega_{k-1}')$.

\bigskip

In the present work, we will study the collection of disks $\cD_n(\omega_1'\cdots \omega_n')$ where the angles $\omega_k^{j_k}$ (for $j_k=1,\dots,4^{k-1}$ and all $k=1,2,\dots,n$) of the tree $\omega_1'\cdots \omega_n'$ are chosen randomly with uniform and independent distributions on the interval $[0,\frac{\pi}{2}]$.

Let us describe this picture once more before moving on further.
The set $\cD_n(\omega_1'\cdots \omega_n')$ consists of $4^n$ disks of radius $4^{-n}$. These can be separated into $4^{n-1}$ groups of four, which are copies of
\[\Omega_n^{j_n}=\bigcup_{\alpha=0}^3\frac{1}{4^{n-1}}T_\alpha^{\omega_n^{j_n}}(\cD_0)\]
translated appropriately within the unit disk so that the ``centre'' of $\Omega_n^{j_n}$ is placed at the centre of the $j_n$-th disk of $\cD_{n-1}(\omega_1'\dots \omega_{n-1}')$ for $j_n=1,2,\dots,4^{n-1}$. That is, each of the $4^{n-1}$ many disks at depth $n-1$ of radius $4^{-(n-1)}$ is replaced by four smaller ones of  radius $4^{-n}$. \Cref{fig:disks} depicts $\cD_n(\omega_1'\cdots \omega_n')$ for $n=0,1,2$.

\section{Favard Length}	\label{sec:Favard_length}

Recall the \emph{Favard length} of a planar set $E\subset\C$ is the integral
\[\Fav(E)=\frac{1}{\pi}\int_0^\pi\left|\proj_\theta E\right|d\theta\]
where $\proj_\theta E$ is the projection of $E$ onto the line with slope $\tan \theta$ passing through the origin, and $|A|$ is the (1-dimensional) Lebesgue measure of $A$.

Now, consider an infinite tree of angles from $[0,\frac{\pi}{2}]$ with root $\omega_1$ and four branches at each vertex, and let $\cD$ be the limit set
\[\cD=\bigcap_{n=0}^\infty\cD_n(\omega_1'\cdots \omega_n').\]
Notice that by construction, $\cD$ a purely unrectifiable planar set. As such, $\Fav(\cD)=0$ and by dominated convergence $\Fav(\cD_n(\omega_1'\cdots \omega_n'))\to0$ while $n\to\infty$. In fact, if the angles are randomly chosen uniformly and independently over $[0,\frac{\pi}{2}]$, by dominated convergence and Fubini $\E[\Fav(\cD)]=0$ and $\E[\Fav(\cD_n(\omega_1'\cdots \omega_n'))]\to0$ as $n\to\infty$, where the expectation is taken over all such angles.

The question arises as to the rate with which $\E[\Fav(\cD_n(\omega_1'\cdots \omega_n'))]$ goes to $0$. This we answer in the following theorem:
\begin{thm} \label{thm:Favard}
	Let $n\in\N$ and consider a tree of angles of height $n$ with each vertex having four branches. Suppose that the angles $\omega_k^{j_k}$ (for all $j_k=1,2,\dots,4^{k-1}$ and all $k=1,2,\dots,n$) are chosen randomly with uniform and independent distributions on the interval $[0,\frac{\pi}{2}]$. Also set $\omega_k'=(\omega_k^1,\omega_k^2,\dots,\omega_k^{4^{k-1}})$ for each $k=1,2,\dots,n$. Then, there exists a constant $c>0$ such that for any $\theta\in[0,\frac{\pi}{2}]$ it holds that
	\begin{equation} \label{eq:average_projection_length}
		\E_{\omega_1'\cdots\, \omega_n'}\left|\proj_\theta\cD_n(\omega_1'\cdots \omega_n')\right|\leq\frac{c}{n}\qquad\forall n\in\N.
	\end{equation}
	
	Consequently,
	\begin{equation} \label{eq:average_Favard_legth}
		\E_{\omega_1'\cdots\, \omega_n'}[\Fav(\cD_n(\omega_1'\cdots \omega_n'))]\leq\frac{c}{n}\qquad\forall n\in\N
	\end{equation}
	and also
	\begin{equation} \label{eq:liminf}
		\liminf_{n\to\infty}n\Fav(\cD_n(\omega_1'\cdots \omega_n'))<\infty\qquad\forall n\in\N\text{ almost surely}.
	\end{equation}
\end{thm}
\noindent Clearly, \eqref{eq:liminf} follows from \eqref{eq:average_Favard_legth} by an immediate application of Fatou's lemma, whereas \eqref{eq:average_Favard_legth} follows from \eqref{eq:average_projection_length} through Fubini.


\section{Statement and use of the main lemma} \label{sec:the_main_lemma}

The present and the following sections are dedicated to the proof of \eqref{eq:average_projection_length}. Towards this goal, we need to introduce \Cref{lem:square_overlap} below, which describes the decay of the average projection when transitioning from depth $k$ to depth $k+1$. The main difficulty will come from obtaining the square factor appearing in \eqref{eq:square_overlap}, which emanates from the naturally occurring overlap of the projections.

\bigskip

From now on, suppose we are given a tree of angles of height $n$ with four branches at each vertex where the angles are uniformly and independently distributed random variables on the interval $[0,\frac{\pi}{2}]$. Recall that given such a tree $\bar{\omega}_{n-k+1}^{j_{n-k+1}}$ is the subtree of height $k$ with the vertex $\omega_{n-k+1}^{j_{n-k+1}}$ as its root. Observe that $\bar{\omega}_1^{j_1}=\omega_1'\cdots \omega_n'$ is the full tree whilst $\bar{\omega}_n^{j_n}=\omega_n^{j_n}$ ($j_n=1,2,\dots,4^{n-1}$) are its leaves, i.e. trees of height $1$.

For any $\theta\in[0,\frac{\pi}{2}]$ and all $k=1,2,\dots,n$, define the following quantities
\begin{align*}
	D_1^{j_n}		&	=\E_{\bar{\omega}_n^{j_n}}\left|\proj_\theta\cD_1(\bar{\omega}_n^{j_n})\right|,	&&	j_n=1,2,\dots,4^{n-1}\\
	D_k^{j_{n-k+1}}	&	=\E_{\bar{\omega}_{n-k+1}^{j_{n-k+1}}}\left|\proj_\theta\cD_k(\bar{\omega}_{n-k+1}^{j_{n-k+1}})\right|,	&&	j_{n-k+1}=1,2,\dots,4^{n-k}\\
	D_n^{j_1}=D_n^1	&	=\E_{\bar{\omega}_1^{j_1}}\left|\proj_\theta\cD_n(\bar{\omega}_1^{j_1})\right|,	&&	j_1=1.
\end{align*}
Notice that, because we are averaging over the independent and identically distributed $\omega_k^{j_k}$,
\[D_k^1=D_k^2=\dots=D_k^{4^{n-k}}\afterline{for any}k=1,2,\dots,n.\]
Therefore, it suffices to work with $D_k^1$; the rest should be identical. Also, note that
\[D_n^1=\E_{\omega_1'\cdots\, \omega_n'}\left|\proj_\theta\cD_n(\omega_1'\cdots \omega_n')\right|.\]

Now, we are ready to state a simple but important lemma. Also, see \cite[Lemma 2.1]{Zha2019}.
\begin{lem} \label{lem:square_overlap}
	With notation as above, if $\omega_k^{j_k}$ are uniformly and independently distributed random variables on $[0,\frac{\pi}{2}]$, there exists a constant $c\geq 4$ such that for any $n\in\N$ (and any $\theta\in[0,\frac{\pi}{2}]$)
	\begin{equation} \label{eq:square_overlap}
		D_{k+1}^1\leq D_k^1-c^{-1}(D_k^1)^2\afterline{for all}k=1,\dots,n-1.
	\end{equation}
\end{lem}
\noindent In our computations later, we will have that $c=12\sqrt2$. But this is possibly not sharp.

Provided this holds true we can give a very compact proof of \Cref{thm:Favard} using induction:
\begin{proof}[Proof of \Cref{thm:Favard}]
	Let $c$ be as in \Cref{lem:square_overlap} and note that $D_2^1\leq D_1^1<2\leq\frac{c}{2}$. Also, $D_1^1<c$.
	
	Next, assume $D_k^1<\frac{c}{k}$ for some $2\leq k\leq n-1$. From \Cref{lem:square_overlap}, and by the monotonicity of the function $x-x^2/c$ in $[0,\frac{c}{2}]$, we see that
	\[D_{k+1}^1\leq D_k^1-c^{-1}(D_k^1)^2<\frac{c}{k}-\frac{c}{k^2}=c\frac{k-1}{k^2}<\frac{c}{k+1}.\]
	Therefore, $D_k^1<\frac{c}{k}$ holds for all for $1\leq k\leq n-1$ and thus for $k=n$ we get
	\[\E_{\omega_1'\cdots\, \omega_n'}\left|\proj_\theta\cD_n(\omega_1'\cdots \omega_n')\right|=D_n^1<\frac{c}{n}.\]
	This is \eqref{eq:average_projection_length}. Equation \eqref{eq:average_Favard_legth} follows after integrating with respect to $\theta$, and \eqref{eq:liminf} after applying Fatou's Lemma.
\end{proof}

\section{Proving the main lemma} \label{sec:proof_of_lemma}

Whatever follows is dedicated to the proof of \eqref{eq:square_overlap}.

\medskip

First, we rewrite the length of the projection of a set in more convenient way. Let $l_\theta$ and $l_\theta^\perp$ be two lines through the origin so that $l_\theta$ forms an angle $\theta$ with the horizontal axis and $l_\theta^\perp$ is perpendicular to $l_\theta$. Also, let $\textbf{n}$ be the unit normal vector of $l_\theta^\perp$. The length of the projection of a planar set $E\subset\C$ onto the line $l_\theta$ can be written as
\begin{equation} \label{eq:projection}
	\left|\proj_\theta E\right|=\left|\{t\in\R\such(l_\theta^\perp+t\textbf{n})\cap E\not=\emptyset\}\right|=\int_{(l_\theta^\perp+t\textbf{n})\cap E\not=\emptyset}dt.
\end{equation}

For brevity, we denote the line $l_\theta^\perp+t\textbf{n}$ by $\needle_\theta(t)$ where $t\in\R$. Additionally, because of the symmetry of our considerations, we can assume without loss of generality that $\theta=0$ ---as we will average over all $\theta$ at the end. So, we can simply omit writing $\theta$ altogether from now on.

\medskip

The idea behind \Cref{lem:square_overlap} is to look at the collection $\cD_n(\omega_1'\cdots \omega_n')$ at depth $n$ but \enquote{zoomed in} so that it looks like depth $1$. Then, we go one level up and look at the disks of depth $n-1$ and $n$ zooming in enough so that they to look like depth $2$; and so forth. If we rewrite the projections in the form of \eqref{eq:projection}, the average overlap at each level is of at least a square factor compared to the total average projection of the level above.

This last comparison is paramount to the proof. It will follow from the fact that the disks in our constructions never get too close to one another. In fact, this observation is not true in the case of the random square Cantor sets, which is the reason why we cannot directly apply the arguments here to the setting of \cite{PerSol2002}.

Let us proceed with the proof of \eqref{eq:square_overlap}.

\bigskip

Fix some $k=1,2,\dots,n$ and recall that by construction (eq. \eqref{eq:first_iteration})
\[\cD_1(\omega_{n-k+1}^1)=\bigcup_{\alpha=0}^3T_\alpha^{\omega_{n-k+1}^1}(\cD_0).\]
This means that each disk from the collection $\cD_k(\bar{\omega}_{n-k+1}^1)$ lies inside one of the above four disks, and therefore we can separate $\cD_k(\bar{\omega}_{n-k+1}^1)$ into four groups of disks depending on their positioning at depth $1$.

\begin{figure}[h]
	\centering
	\begin{tikzpicture}[scale=4,trim right=(RRR),trim left=(LLL)]
		\node (RRR) at (1,0) {};
		\node (LLL) at (-1,0) {};
		\draw [opacity=0.1] (0,0) circle (1);
		\node (DDD) at (-1.35,-0.15) {$\cD_k(\bar{\omega}_{n-k+1}^1)$};
		
		\fill [opacity=0.05] ({3/4*cos(60)}, {3/4*sin(60)})	circle (1/4);
		\fill [opacity=0.05] ({3/4*cos(150)},{3/4*sin(150)}) circle (1/4);
		\fill [opacity=0.05] ({3/4*cos(240)},{3/4*sin(240)}) circle (1/4);
		\fill [opacity=0.05] ({3/4*cos(330)},{3/4*sin(330)}) circle (1/4);
		\draw [opacity=0.4] (0,0) -- (0,1/3.5);
		\draw [opacity=0.4] (0,1/4) arc (90:60:1/4) node [opacity=1,right] {${\omega}_{n-k+1}^1$};
		\draw [opacity=0.4] ({1/2*cos(60)},	{1/2*sin(60)}) -- ({1/2*cos(240)},{1/2*sin(240)});
		\draw [opacity=0.4] ({1/2*cos(150)},{1/2*sin(150)}) -- ({1/2*cos(-30)},{1/2*sin(-30)});
		\draw [opacity=0.4] ({3/4*cos(60)+1/8*cos(10)},{3/4*sin(60)+1/8*sin(10)}) -- ({3/4*cos(60)},{3/4*sin(60)}) -- ({3/4*cos(60)},{3/4*sin(60)+1/12});
		\draw [opacity=0.4] ({3/4*cos(60)},{3/4*sin(60)+1/16}) arc (90:10:1/16) node [midway] (AAA) {};
		\node at ({3/4*cos(60)},{3/4*sin(60)}) (BBB) {};
		\draw [opacity=0] (AAA) -- (BBB) node [below=-3pt,right=-3pt] (CCC) {};
		\draw [Stealth-,very thin] (CCC) -- (0.6,0.4) node [right] {$\omega_{n-k+2}^2$};
		
		\fill [opacity=0.1] ({3/16*cos(30)}, {3/16*sin(30)}) +({3/4*cos(-30)},{3/4*sin(-30)}) circle (1/4^2);
		\fill [opacity=0.1] ({3/16*cos(120)},{3/16*sin(120)})+({3/4*cos(-30)},{3/4*sin(-30)}) circle (1/4^2);
		\fill [opacity=0.1] ({3/16*cos(210)},{3/16*sin(210)})+({3/4*cos(-30)},{3/4*sin(-30)}) circle (1/4^2);
		\fill [opacity=0.1] ({3/16*cos(300)},{3/16*sin(300)})+({3/4*cos(-30)},{3/4*sin(-30)}) circle (1/4^2);
		\fill [opacity=0.1] ({3/16*cos(10)}, {3/16*sin(10)}) +({3/4*cos(60)},{3/4*sin(60)}) circle (1/4^2);
		\fill [opacity=0.1] ({3/16*cos(100)},{3/16*sin(100)})+({3/4*cos(60)},{3/4*sin(60)}) circle (1/4^2);
		\fill [opacity=0.1] ({3/16*cos(190)},{3/16*sin(190)})+({3/4*cos(60)},{3/4*sin(60)}) circle (1/4^2);
		\fill [opacity=0.1] ({3/16*cos(280)},{3/16*sin(280)})+({3/4*cos(60)},{3/4*sin(60)}) circle (1/4^2);
		\fill [opacity=0.1] ({3/64*cos(20)},{3/64*sin(20)})+({3/16*cos(10)+3/4*cos(60)},{3/16*sin(10)+3/4*sin(60)}) circle (1/4^3);
		\fill [opacity=0.1] ({3/64*cos(110)},{3/64*sin(110)})+({3/16*cos(10)+3/4*cos(60)},{3/16*sin(10)+3/4*sin(60)}) circle (1/4^3);
		\fill [opacity=0.1] ({3/64*cos(200)},{3/64*sin(200)})+({3/16*cos(10)+3/4*cos(60)},{3/16*sin(10)+3/4*sin(60)}) circle (1/4^3);
		\fill [opacity=0.1] ({3/64*cos(290)},{3/64*sin(290)})+({3/16*cos(10)+3/4*cos(60)},{3/16*sin(10)+3/4*sin(60)}) circle (1/4^3);
		\fill [opacity=0.1] ({3/64*cos(-5)},{3/64*sin(-5)})+({3/16*cos(100)+3/4*cos(60)},{3/16*sin(100)+3/4*sin(60)}) circle (1/4^3);
		\fill [opacity=0.1] ({3/64*cos(85)},{3/64*sin(85)})+({3/16*cos(100)+3/4*cos(60)},{3/16*sin(100)+3/4*sin(60)}) circle (1/4^3);
		\fill [opacity=0.1] ({3/64*cos(175)},{3/64*sin(175)})+({3/16*cos(100)+3/4*cos(60)},{3/16*sin(100)+3/4*sin(60)}) circle (1/4^3);
		\fill [opacity=0.1] ({3/64*cos(265)},{3/64*sin(265)})+({3/16*cos(100)+3/4*cos(60)},{3/16*sin(100)+3/4*sin(60)}) circle (1/4^3);
		\fill [opacity=0.1] ({3/64*cos(45)},{3/64*sin(45)})+({3/16*cos(190)+3/4*cos(60)},{3/16*sin(190)+3/4*sin(60)}) circle (1/4^3);
		\fill [opacity=0.1] ({3/64*cos(135)},{3/64*sin(135)})+({3/16*cos(190)+3/4*cos(60)},{3/16*sin(190)+3/4*sin(60)}) circle (1/4^3);
		\fill [opacity=0.1] ({3/64*cos(225)},{3/64*sin(225)})+({3/16*cos(190)+3/4*cos(60)},{3/16*sin(190)+3/4*sin(60)}) circle (1/4^3);
		\fill [opacity=0.1] ({3/64*cos(315)},{3/64*sin(315)})+({3/16*cos(190)+3/4*cos(60)},{3/16*sin(190)+3/4*sin(60)}) circle (1/4^3);
		\fill [opacity=0.1] ({3/64*cos(65)},{3/64*sin(65)})+({3/16*cos(280)+3/4*cos(60)},{3/16*sin(280)+3/4*sin(60)}) circle (1/4^3);
		\fill [opacity=0.1] ({3/64*cos(155)},{3/64*sin(155)})+({3/16*cos(280)+3/4*cos(60)},{3/16*sin(280)+3/4*sin(60)}) circle (1/4^3);
		\fill [opacity=0.1] ({3/64*cos(245)},{3/64*sin(245)})+({3/16*cos(280)+3/4*cos(60)},{3/16*sin(280)+3/4*sin(60)}) circle (1/4^3);
		\fill [opacity=0.1] ({3/64*cos(335)},{3/64*sin(335)})+({3/16*cos(280)+3/4*cos(60)},{3/16*sin(280)+3/4*sin(60)}) circle (1/4^3);
		\fill [opacity=0.1] ({3/16*cos(45)}, {3/16*sin(45)}) +({3/4*cos(150)},{3/4*sin(150)}) circle (1/4^2);
		\fill [opacity=0.1] ({3/16*cos(135)},{3/16*sin(135)})+({3/4*cos(150)},{3/4*sin(150)}) circle (1/4^2);
		\fill [opacity=0.1] ({3/16*cos(225)},{3/16*sin(225)})+({3/4*cos(150)},{3/4*sin(150)}) circle (1/4^2);
		\fill [opacity=0.1] ({3/16*cos(315)},{3/16*sin(315)})+({3/4*cos(150)},{3/4*sin(150)}) circle (1/4^2);
		\fill [opacity=0.1] ({3/16*cos(60)}, {3/16*sin(60)}) +({3/4*cos(240)},{3/4*sin(240)}) circle (1/4^2);
		\fill [opacity=0.1] ({3/16*cos(150)},{3/16*sin(150)})+({3/4*cos(240)},{3/4*sin(240)}) circle (1/4^2);
		\fill [opacity=0.1] ({3/16*cos(240)},{3/16*sin(240)})+({3/4*cos(240)},{3/4*sin(240)}) circle (1/4^2);
		\fill [opacity=0.1] ({3/16*cos(330)},{3/16*sin(330)})+({3/4*cos(240)},{3/4*sin(240)}) circle (1/4^2);
		
		\draw[ultra thin] ({3/4*cos(60)},	{3/4*sin(60)})	circle (1/4);
		\draw[ultra thin] ({3/4*cos(150)},	{3/4*sin(150)}) circle (1/4);
		\draw[ultra thin] ({3/4*cos(240)},	{3/4*sin(240)}) circle (1/4);
		\draw[ultra thin] ({3/4*cos(330)},	{3/4*sin(330)}) circle (1/4);
		
		\node at ({5/4*cos(330)},{5/4*sin(330)}) {$\cT_0^k(\bar{\omega}_{n-k+1}^1)$};
		\node at ({5/4*cos(60)+1/2},{5/4*sin(60)-1/12}) {$\cT_1^k(\bar{\omega}_{n-k+1}^1)=\cT_1^k({\omega}_{n-k+1}^1,\bar{\omega}_{n-k+2}^2)$};
		\node at ({5/4*cos(150)},{5/4*sin(150)})  {$\cT_2^k(\bar{\omega}_{n-k+1}^1)$};
		\node at ({5/4*cos(240)-1/4},{5/4*sin(240)+1/4}) {$\cT_3^k(\bar{\omega}_{n-k+1}^1)$};
	\end{tikzpicture}
	\caption{The four groups of disks at depth $k$ rotated by $\omega_{n-k+1}^1$, which make up~$\cD_k(\bar{\omega}_{n-k+1}^1)$.}
	\label{fig:grouping_compass}
\end{figure}
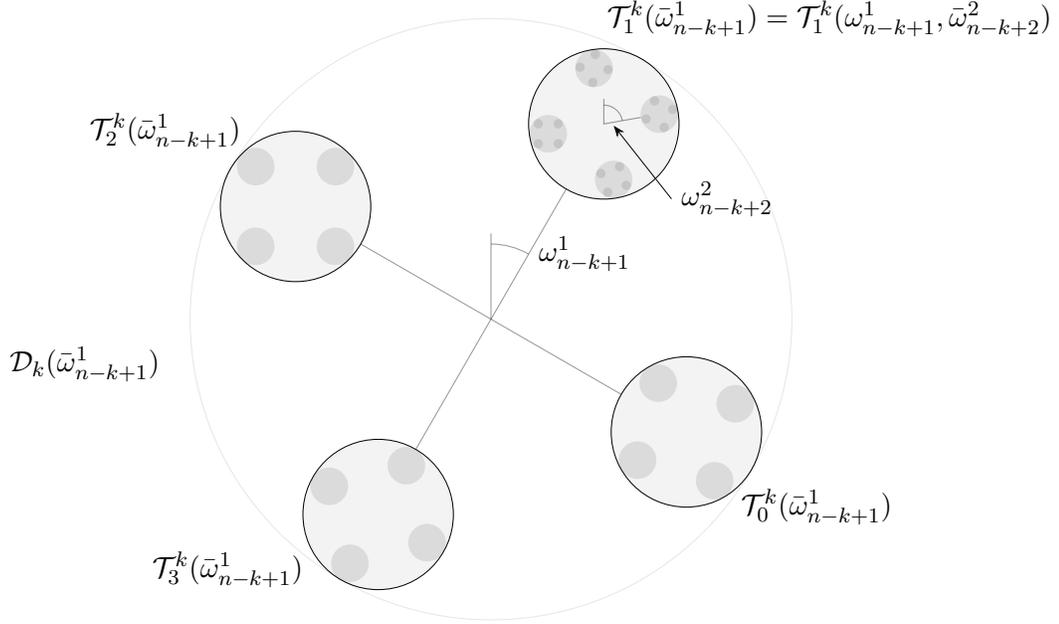

More precisely, for each $\alpha=0,1,2,3$ define $\cT_\alpha^k(\bar{\omega}_{n-k+1}^1)$ as
\begin{equation}	\label{def:compass}
	\cT_\alpha^k(\bar{\omega}_{n-k+1}^1)=T_\alpha^{\omega_{n-k+1}^1}(\cD_0)\bigcap\cD_k(\bar{\omega}_{n-k+1}^1).
\end{equation}
That is, the set $\cT_\alpha^k(\bar{\omega}_{n-k+1}^1)$ consists of those disks of $\cD_k(\bar{\omega}_{n-k+1}^1)$ which lie inside the $\frac{1}{4}$-radius disk $T_\alpha^{\omega_{n-k+1}^1}(\cD_0)$ as in \Cref{fig:grouping_compass}. We can think of $\cT_\alpha^k(\bar{\omega}_{n-k+1}^1)$ as the East, North, West, and South parts of $\cD_k(\bar{\omega}_{n-k+1}^1)$, respectively for $\alpha=0,1,2,3$. From this definition, it is also clear that
\begin{equation} \label{eq:separate_the_disks}
	\cD_k(\bar{\omega}_{n-k+1}^1)=\bigcup_{\alpha=0}^3\cT_\alpha^k(\bar{\omega}_{n-k+1}^1).
\end{equation}

In fact, $\cT_\alpha^k(\bar{\omega}_{n-k+1}^1)$ depends only on the angle $\omega_{n-k+1}^1$ and the subtree $\bar{\omega}_{n-k+2}^{4\cdot1-3+\alpha}=\bar{\omega}_{n-k+2}^{1+\alpha}$. (Recall our enumeration of the angles in \Cref{sec:the_disks}.) Thus, we can write $\cT_\alpha^k(\bar{\omega}_{n-k+1}^1)$ as
\begin{equation}	\label{eq:dependence}
	\cT_\alpha^k(\bar{\omega}_{n-k+1}^1)=\cT_\alpha^k(\omega_{n-k+1}^1,\bar{\omega}_{n-k+2}^{1+\alpha}).
\end{equation}

\subsection{Key observations} \label{subsec:key}
There are two key observations regarding the sets $\cT_\alpha^k(\bar{\omega}_{n-k+1}^1)$. First, note that each point of the interval $(-1,1)$ can be covered by at most two of the projections $\proj T_\alpha^{\omega_{n-k+1}^1}(\cD_0)$ for different $\alpha$'s. Since $\cT_\alpha^k(\bar{\omega}_{n-k+1}^1)\subset T_\alpha^{\omega_{n-k+1}^1}(\cD_0)$, the same holds true for $\proj\cT_\alpha^k(\bar{\omega}_{n-k+1}^1)$; the intersection $\bigcap_\alpha\proj\cT_\alpha^k(\bar{\omega}_{n-k+1}^1)$ is empty when the intersection is over more than two values of $\alpha$.

Second, we can compare the average projections of $\cT_\alpha^k(\bar{\omega}_{n-k+1}^1)=\cT_\alpha^k(\omega_{n-k+1}^1,\bar{\omega}_{n-k+2}^{1+\alpha})$ and $\cD_{k-1}(\bar{\omega}_{n-k+2}^{1+\alpha})$. Notice that both these collections consist of $4^{k-1}$ many disks, which in fact have the same $n$-depth enumerations. This means that they correspond to the same disks of the collection $\cD_n(\omega_1'\cdots \omega_n')$. The difference is that the disks of the former are translated according to $\cD_1(\omega_{n-k+1}^1)$ and have radius $4^{-k}$, whereas the ones of the latter have radius $4^{-(k-1)}$.

\begin{figure}[h]
	\centering
	\begin{tikzpicture}[scale=2]
		\draw [opacity=0.1] (0,0) circle (1);
		\node at (-1.3,-0.2) {$\cD_k(\bar{\omega}_{n-k+1}^1)$};
		
		\fill [opacity=0.05] ({3/4*cos(60)}, {3/4*sin(60)})	circle (1/4);
		\fill [opacity=0.05] ({3/4*cos(150)},{3/4*sin(150)}) circle (1/4);
		\fill [opacity=0.05] ({3/4*cos(240)},{3/4*sin(240)}) circle (1/4);
		\fill [opacity=0.05] ({3/4*cos(330)},{3/4*sin(330)}) circle (1/4);
		\draw [opacity=0.4] (0,0) -- (0,1/3.5) (0,1/4) arc (90:60:1/4) node [opacity=1,right] {\scriptsize ${\omega}_{n-k+1}^1$};
		\draw [opacity=0.4] ({1/2*cos(60)},	{1/2*sin(60)}) -- (0,0);
		
		\fill [opacity=0.1] ({3/16*cos(30)}, {3/16*sin(30)}) +({3/4*cos(-30)},{3/4*sin(-30)}) circle (1/4^2);
		\fill [opacity=0.1] ({3/16*cos(120)},{3/16*sin(120)})+({3/4*cos(-30)},{3/4*sin(-30)}) circle (1/4^2);
		\fill [opacity=0.1] ({3/16*cos(210)},{3/16*sin(210)})+({3/4*cos(-30)},{3/4*sin(-30)}) circle (1/4^2);
		\fill [opacity=0.1] ({3/16*cos(300)},{3/16*sin(300)})+({3/4*cos(-30)},{3/4*sin(-30)}) circle (1/4^2);
		\fill [opacity=0.1] ({3/16*cos(10)}, {3/16*sin(10)}) +({3/4*cos(60)},{3/4*sin(60)}) circle (1/4^2);
		\fill [opacity=0.1] ({3/16*cos(100)},{3/16*sin(100)})+({3/4*cos(60)},{3/4*sin(60)}) circle (1/4^2);
		\fill [opacity=0.1] ({3/16*cos(190)},{3/16*sin(190)})+({3/4*cos(60)},{3/4*sin(60)}) circle (1/4^2);
		\fill [opacity=0.1] ({3/16*cos(280)},{3/16*sin(280)})+({3/4*cos(60)},{3/4*sin(60)}) circle (1/4^2);
		\fill [opacity=0.1] ({3/64*cos(20)},{3/64*sin(20)})+({3/16*cos(10)+3/4*cos(60)},{3/16*sin(10)+3/4*sin(60)}) circle (1/4^3);
		\fill [opacity=0.1] ({3/64*cos(110)},{3/64*sin(110)})+({3/16*cos(10)+3/4*cos(60)},{3/16*sin(10)+3/4*sin(60)}) circle (1/4^3);
		\fill [opacity=0.1] ({3/64*cos(200)},{3/64*sin(200)})+({3/16*cos(10)+3/4*cos(60)},{3/16*sin(10)+3/4*sin(60)}) circle (1/4^3);
		\fill [opacity=0.1] ({3/64*cos(290)},{3/64*sin(290)})+({3/16*cos(10)+3/4*cos(60)},{3/16*sin(10)+3/4*sin(60)}) circle (1/4^3);
		\fill [opacity=0.1] ({3/64*cos(-5)},{3/64*sin(-5)})+({3/16*cos(100)+3/4*cos(60)},{3/16*sin(100)+3/4*sin(60)}) circle (1/4^3);
		\fill [opacity=0.1] ({3/64*cos(85)},{3/64*sin(85)})+({3/16*cos(100)+3/4*cos(60)},{3/16*sin(100)+3/4*sin(60)}) circle (1/4^3);
		\fill [opacity=0.1] ({3/64*cos(175)},{3/64*sin(175)})+({3/16*cos(100)+3/4*cos(60)},{3/16*sin(100)+3/4*sin(60)}) circle (1/4^3);
		\fill [opacity=0.1] ({3/64*cos(265)},{3/64*sin(265)})+({3/16*cos(100)+3/4*cos(60)},{3/16*sin(100)+3/4*sin(60)}) circle (1/4^3);
		\fill [opacity=0.1] ({3/64*cos(45)},{3/64*sin(45)})+({3/16*cos(190)+3/4*cos(60)},{3/16*sin(190)+3/4*sin(60)}) circle (1/4^3);
		\fill [opacity=0.1] ({3/64*cos(135)},{3/64*sin(135)})+({3/16*cos(190)+3/4*cos(60)},{3/16*sin(190)+3/4*sin(60)}) circle (1/4^3);
		\fill [opacity=0.1] ({3/64*cos(225)},{3/64*sin(225)})+({3/16*cos(190)+3/4*cos(60)},{3/16*sin(190)+3/4*sin(60)}) circle (1/4^3);
		\fill [opacity=0.1] ({3/64*cos(315)},{3/64*sin(315)})+({3/16*cos(190)+3/4*cos(60)},{3/16*sin(190)+3/4*sin(60)}) circle (1/4^3);
		\fill [opacity=0.1] ({3/64*cos(65)},{3/64*sin(65)})+({3/16*cos(280)+3/4*cos(60)},{3/16*sin(280)+3/4*sin(60)}) circle (1/4^3);
		\fill [opacity=0.1] ({3/64*cos(155)},{3/64*sin(155)})+({3/16*cos(280)+3/4*cos(60)},{3/16*sin(280)+3/4*sin(60)}) circle (1/4^3);
		\fill [opacity=0.1] ({3/64*cos(245)},{3/64*sin(245)})+({3/16*cos(280)+3/4*cos(60)},{3/16*sin(280)+3/4*sin(60)}) circle (1/4^3);
		\fill [opacity=0.1] ({3/64*cos(335)},{3/64*sin(335)})+({3/16*cos(280)+3/4*cos(60)},{3/16*sin(280)+3/4*sin(60)}) circle (1/4^3);
		\fill [opacity=0.1] ({3/16*cos(45)}, {3/16*sin(45)}) +({3/4*cos(150)},{3/4*sin(150)}) circle (1/4^2);
		\fill [opacity=0.1] ({3/16*cos(135)},{3/16*sin(135)})+({3/4*cos(150)},{3/4*sin(150)}) circle (1/4^2);
		\fill [opacity=0.1] ({3/16*cos(225)},{3/16*sin(225)})+({3/4*cos(150)},{3/4*sin(150)}) circle (1/4^2);
		\fill [opacity=0.1] ({3/16*cos(315)},{3/16*sin(315)})+({3/4*cos(150)},{3/4*sin(150)}) circle (1/4^2);
		\fill [opacity=0.1] ({3/16*cos(60)}, {3/16*sin(60)}) +({3/4*cos(240)},{3/4*sin(240)}) circle (1/4^2);
		\fill [opacity=0.1] ({3/16*cos(150)},{3/16*sin(150)})+({3/4*cos(240)},{3/4*sin(240)}) circle (1/4^2);
		\fill [opacity=0.1] ({3/16*cos(240)},{3/16*sin(240)})+({3/4*cos(240)},{3/4*sin(240)}) circle (1/4^2);
		\fill [opacity=0.1] ({3/16*cos(330)},{3/16*sin(330)})+({3/4*cos(240)},{3/4*sin(240)}) circle (1/4^2);
		
		\draw[ultra thin] ({3/4*cos(60)},	{3/4*sin(60)})	circle (1/4);
		\draw[ultra thin] ({3/4*cos(150)},	{3/4*sin(150)}) circle (1/4);
		\draw[ultra thin] ({3/4*cos(240)},	{3/4*sin(240)}) circle (1/4);
		\draw[ultra thin] ({3/4*cos(330)},	{3/4*sin(330)}) circle (1/4);
		
		\node at ({5/4*cos(330)-1/12},	{6/4*sin(330)-1/12}){\scriptsize $\cT_0^k(\bar{\omega}_{n-k+1}^1)$};
		\node at ({5/4*cos(60)+1/4},	{5/4*sin(60)-1/12})	{$\cT_1^k(\bar{\omega}_{n-k+1}^1)$};
		\node at ({5/4*cos(150)+1/6},	{5/4*sin(150)+1/8})	{\scriptsize $\cT_2^k(\bar{\omega}_{n-k+1}^1)$};
		\node at ({5/4*cos(240)-1/16},	{5/4*sin(240)+1/16}){\scriptsize $\cT_3^k(\bar{\omega}_{n-k+1}^1)$};
		
		\draw[ultra thin] (4,0) circle (1);
		\node at ({6/4*cos(60)+4+1/4},	{6/4*sin(60)-1/3.5}) {$\cD_{k-1}(\bar{\omega}_{n-k+2}^2)$};
		\draw ({3/4*cos(60)+5/16},{3/4*sin(60)-1/16}) -- ({4-3/4*cos(60)-5/16}, {3/4*sin(60)+4/16});
		\draw ({3/4*cos(60)+5/16},{3/4*sin(60)-1/16}) -- ({4-3/4*cos(60)-5/16}, {-3/4*sin(60)-4/16});
		\draw [opacity=0.4] ({1/2*cos(20)+4},{1/2*sin(20)}) -- (4,0) -- (4,1/3);
		\draw [opacity=0.4] (4,1/4) arc (90:20:1/4) node [opacity=1,right=3pt,below=3pt] {${\omega}_{n-k+2}^2$};
		\fill [opacity=0.05] ({3/4*cos(20)+4}, {3/4*sin(20)})  circle (1/4);
		\fill [opacity=0.05] ({3/4*cos(110)+4},{3/4*sin(110)}) circle (1/4);
		\fill [opacity=0.05] ({3/4*cos(200)+4},{3/4*sin(200)}) circle (1/4);
		\fill [opacity=0.05] ({3/4*cos(290)+4},{3/4*sin(290)}) circle (1/4);
		\fill [opacity=0.1] ({3/16*cos(20)+4},{3/16*sin(20)})+({3/4*cos(20)},{3/4*sin(20)}) circle (1/4^2);
		\fill [opacity=0.1] ({3/16*cos(110)+4},{3/16*sin(110)})+({3/4*cos(20)},{3/4*sin(20)}) circle (1/4^2);
		\fill [opacity=0.1] ({3/16*cos(200)+4},{3/16*sin(200)})+({3/4*cos(20)},{3/4*sin(20)}) circle (1/4^2);
		\fill [opacity=0.1] ({3/16*cos(290)+4},{3/16*sin(290)})+({3/4*cos(20)},{3/4*sin(20)}) circle (1/4^2);
		\fill [opacity=0.1] ({3/16*cos(-5)+4},{3/16*sin(-5)})+({3/4*cos(110)},{3/4*sin(110)}) circle (1/4^2);
		\fill [opacity=0.1] ({3/16*cos(85)+4},{3/16*sin(85)})+({3/4*cos(110)},{3/4*sin(110)}) circle (1/4^2);
		\fill [opacity=0.1] ({3/16*cos(175)+4},{3/16*sin(175)})+({3/4*cos(110)},{3/4*sin(110)}) circle (1/4^2);
		\fill [opacity=0.1] ({3/16*cos(265)+4},{3/16*sin(265)})+({3/4*cos(110)},{3/4*sin(110)}) circle (1/4^2);
		\fill [opacity=0.1] ({3/16*cos(45)+4},{3/16*sin(45)})+({3/4*cos(200)},{3/4*sin(200)}) circle (1/4^2);
		\fill [opacity=0.1] ({3/16*cos(135)+4},{3/16*sin(135)})+({3/4*cos(200)},{3/4*sin(200)}) circle (1/4^2);
		\fill [opacity=0.1] ({3/16*cos(225)+4},{3/16*sin(225)})+({3/4*cos(200)},{3/4*sin(200)}) circle (1/4^2);
		\fill [opacity=0.1] ({3/16*cos(315)+4},{3/16*sin(315)})+({3/4*cos(200)},{3/4*sin(200)}) circle (1/4^2);
		\fill [opacity=0.1] ({3/16*cos(65)+4},{3/16*sin(65)})+({3/4*cos(290)},{3/4*sin(290)}) circle (1/4^2);
		\fill [opacity=0.1] ({3/16*cos(155)+4},{3/16*sin(155)})+({3/4*cos(290)},{3/4*sin(290)}) circle (1/4^2);
		\fill [opacity=0.1] ({3/16*cos(245)+4},{3/16*sin(245)})+({3/4*cos(290)},{3/4*sin(290)}) circle (1/4^2);
		\fill [opacity=0.1] ({3/16*cos(335)+4},{3/16*sin(335)})+({3/4*cos(290)},{3/4*sin(290)}) circle (1/4^2);
	\end{tikzpicture}
	\caption{Dilating $\cT_1^k(\bar{\omega}_{n-k+1})$ by 4 gives a copy of $\cD_{k-1}(\bar{\omega}_{n-k+2}^2)$.}
	\label{fig:zooming_in}
\end{figure}
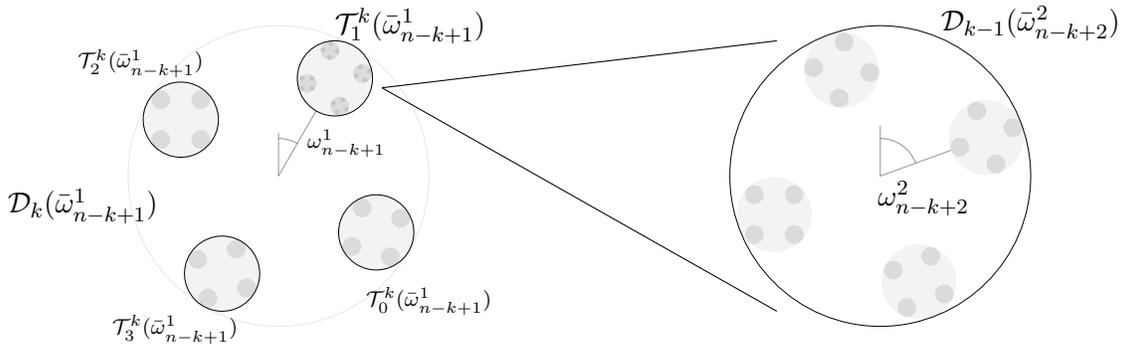

Consequently, $\cT_\alpha^k(\bar{\omega}_{n-k+1}^1)$ is a shifted copy of $\cD_{k-1}(\bar{\omega}_{n-k+2}^{1+\alpha})$ dilated by a factor of $1/4$. (See \Cref{fig:zooming_in}.) As such, the (average of the) projections of $\cT_\alpha^k(\bar{\omega}_{n-k+1}^1)$ and $\cD_{k-1}(\bar{\omega}_{n-k+2}^{1+\alpha})$ should also differ by a factor of $1/4$. In other words, for any $\alpha=0,1,2,3$ we have
\begin{equation} \label{eq:zoom}
	\begin{split}
		\E_{\bar{\omega}_{n-k+1}^1}\left|\proj\cT_\alpha^k(\bar{\omega}_{n-k+1}^1)\right|
			&	{}=\E_{\omega_{n-k+1}^1}\E_{\bar{\omega}_{n-k+2}^{1+\alpha}}\left|\proj\cT_\alpha^k(\omega_{n-k+1}^1,\bar{\omega}_{n-k+2}^{1+\alpha})\right|\\
			&	=\frac{1}{4}\E_{\bar{\omega}_{n-k+2}^{1+\alpha}}\left|\proj\cD_{k-1}(\bar{\omega}_{n-k+2}^{1+\alpha})\right|.
	\end{split}
\end{equation}

\subsection{The estimates} \label{subsec:estimates}
Utilising the above, we can now estimate $D_k^1$ in terms of $D_{k-1}^1$.

For starters, note that from \eqref{eq:separate_the_disks} we can write
\begin{align*}
	D_k^1	&	=\E_{\bar{\omega}_{n-k+1}^1}\left|\proj\cD_k(\bar{\omega}_{n-k+1}^1)\right|\\
		\begin{split}
			&	=\E_{\bar{\omega}_{n-k+1}^1}\sum_{\alpha=0}^3\left|\proj\cT_\alpha^k(\bar{\omega}_{n-k+1}^1)\right|	-\E_{\bar{\omega}_{n-k+1}^1}\sum_{\substack{\alpha,\beta=0\\\alpha\not=\beta}}^3\left|\proj\cT_\alpha^k(\bar{\omega}_{n-k+1}^1)\cap\proj\cT_\beta^k(\bar{\omega}_{n-k+1}^1)\right|	+\\
			&\quad+\E_{\bar{\omega}_{n-k+1}^1}\sum_{\substack{\alpha,\beta,\gamma=0\\\alpha\not=\beta\not=\gamma\not=\alpha}}^3\left|\proj\cT_\alpha^k(\bar{\omega}_{n-k+1}^1)\cap\proj\cT_\beta^k(\bar{\omega}_{n-k+1}^1)\cap\proj\cT_\gamma^k(\bar{\omega}_{n-k+1}^1)\right|	-\\
			&\quad-\E_{\bar{\omega}_{n-k+1}^1}\left|\proj\cT_0^k(\bar{\omega}_{n-k+1}^1)\cap\proj\cT_1^k(\bar{\omega}_{n-k+1}^1)\cap\proj\cT_2^k(\bar{\omega}_{n-k+1}^1)\cap\proj\cT_3^k(\bar{\omega}_{n-k+1}^1)\right|.
		\end{split}
\end{align*}
The last two lines equal $0$ from our first observation above (in \Cref{subsec:key}). Furthermore, we can disregard all but one of the summands from the second sum to get an inequality:
\begin{equation}	\label{eq:the_first_inequality}
	D_k^1\leq\E_{\bar{\omega}_{n-k+1}^1}\sum_{\alpha=0}^3\left|\proj\cT_\alpha^k(\bar{\omega}_{n-k+1}^1)\right|
		-\E_{\bar{\omega}_{n-k+1}^1}\left|\proj\cT_0^k(\bar{\omega}_{n-k+1}^1)\cap\proj\cT_1^k(\bar{\omega}_{n-k+1}^1)\right|.
\end{equation}
This last step might seem rather crude, but it will suffice for our purposes. Besides, \Cref{thm:Favard} eventually establishes an equality considering Mattila's lower bound.

Utilising \eqref{eq:zoom}, we see that
\begin{align*}
	\E_{\bar{\omega}_{n-k+1}^1}\sum_{\alpha=0}^3\left|\proj\cT_\alpha^k(\bar{\omega}_{n-k+1}^1)\right|
		&	=\frac{1}{4}\sum_{\alpha=0}^3\E_{\bar{\omega}_{n-k+2}^{1+\alpha}}\left|\proj\cD_{k-1}(\bar{\omega}_{n-k+2}^{1+\alpha})\right|\\
		&	=\frac{1}{4}(D_{k-1}^1+D_{k-1}^2+D_{k-1}^3+D_{k-1}^4)\\
		&	=D_{k-1}^1,
\end{align*}
since $D_{k-1}^{1+\alpha}=D_{k-1}^1$ for any $\alpha=0,1,2,3$. Applying this to \eqref{eq:the_first_inequality}, we get
\begin{equation} \label{eq:extra_overlap}
	D_k^1	\leq D_{k-1}^1	-\E_{\bar{\omega}_{n-k+1}^1}\left|\proj\cT_0^k(\bar{\omega}_{n-k+1}^1)\cap\proj\cT_1^k(\bar{\omega}_{n-k+1}^1)\right|.
\end{equation}

\medskip

The final big step is to estimate the overlap term $\left|\proj\cT_0^k(\bar{\omega}_{n-k+1}^1)\cap\proj\cT_1^k(\bar{\omega}_{n-k+1}^1)\right|$ from below. But recall that $\cT_0^k(\bar{\omega}_{n-k+1}^1)$ and $\cT_1^k(\bar{\omega}_{n-k+1}^1)$ depend (aside from $\omega_{n-k+1}^1$) respectively on $\bar{\omega}_{n-k+2}^1$ and $\bar{\omega}_{n-k+2}^2$ as in \eqref{eq:dependence}.


First, we average with respect to the subtrees $\bar{\omega}_{n-k+2}^1$ and $\bar{\omega}_{n-k+2}^2$, and afterwards we integrate over their common ancestor $\omega_{n-k+1}^1$. To simplify the notation, let us write $\bar{\omega}_{n-k+2}^{1,2}$ for both the subtrees $\bar{\omega}_{n-k+2}^1$ and $\bar{\omega}_{n-k+2}^2$. Then, we have
\begin{align*}
	&	\E_{\bar{\omega}_{n-k+2}^{1,2}}\left|\proj\cT_0^k(\bar{\omega}_{n-k+1}^1)\cap\proj\cT_1^k(\bar{\omega}_{n-k+1}^1)\right|\\
	&	{}=\E_{\bar{\omega}_{n-k+2}^{1,2}}\left|\proj\cT_0^k(\omega_{n-k+1}^1,\bar{\omega}_{n-k+2}^1)\cap\proj\cT_1^k(\omega_{n-k+1}^1,\bar{\omega}_{n-k+2}^2)\right|\\
	&	{}\xlongequal{\eqref{eq:projection}}	\int\P_{\bar{\omega}_{n-k+2}^{1,2}}\left(\needle(t)\cap\cT_0^k(\omega_{n-k+1}^1,\bar{\omega}_{n-k+2}^1)\not=\emptyset	\textsc{ and }	\needle(t)\cap\cT_1^k(\omega_{n-k+1}^1,\bar{\omega}_{n-k+2}^2)\not=\emptyset\right)dt\\
	&	{}=\int\P_{\bar{\omega}_{n-k+2}^{1,2}}\left(\needle(t)\cap\cT_0^k(\omega_{n-k+1}^1,\bar{\omega}_{n-k+2}^1)\not=\emptyset\right)	\!\cdot\!
	\P_{\bar{\omega}_{n-k+2}^{1,2}}\left(\needle(t)\cap\cT_1^k(\omega_{n-k+1}^1,\bar{\omega}_{n-k+2}^2)\not=\emptyset\right)dt\\
	&	{}=\int\P_{\bar{\omega}_{n-k+2}^1}\left(\needle(t)\cap\cT_0^k(\omega_{n-k+1}^1,\bar{\omega}_{n-k+2}^1)\not=\emptyset\right)	\!\cdot\!	
	\P_{\bar{\omega}_{n-k+2}^2}\left(\needle(t)\cap\cT_1^k(\omega_{n-k+1}^1,\bar{\omega}_{n-k+2}^2)\not=\emptyset\right)dt
	\\
	&	{}=:\cE(\omega_{n-k+1}^1).
\end{align*}
The 3rd equality above holds because for a fixed angle $\omega_{n-k+1}^1$
the events
\[\{\needle(t)\cap\cT_\alpha^k(\omega_{n-k+1}^1,\bar{\omega}_{n-k+2}^{1+\alpha})\not=\emptyset\}\]
for different $\alpha$'s are independent.

It would be very nice if these two events would have the same probability for each $t$.
Then, we would use H{\"o}lder's inequality to get that
\begin{align*}
	\cE(\omega_{n-k+1}^1)
	&	=	\int\left[\P_{\bar{\omega}_{n-k+2}^1}\left(\needle(t)\cap\cT_0^k(\omega_{n-k+1}^1,\bar{\omega}_{n-k+2}^1)\not=\emptyset\right)\right]^2dt\\
	&	\geq	C\left(\int\P_{\bar{\omega}_{n-k+2}^1}\left(\needle(t)\cap\cT_0^k(\omega_{n-k+1}^1,\bar{\omega}_{n-k+2}^1)\not=\emptyset\right)dt\right)^2.
\end{align*}
Unfortunately, this is not the case.

\medskip

Nevertheless, there is a way to approximate the overlap $\cE(\omega_{n-k+1}^1)$.

To keep things clean, let us denote
\begin{equation*}
	\psi:=\omega_{n-k+1}^1,\ \bar{\psi}^1:=\bar{\omega}_{n-k+2}^1,\ \bar{\psi}^2:=\bar{\omega}_{n-k+2}^2 \inline{and} s(\psi):=\frac{3}{4}(1-\cos \psi),
\end{equation*}
and keep all angles fixed for now.

As discussed in \Cref{subsec:key}, the average projections of $\cT_0^k(\psi,\bar{\psi}^{1})$ and $\cT_1^k(\psi,\bar{\psi}^{2})$ are shifted (and dilated) copies of $\cD_{k-1}(\bar{\psi}^{1})$ and $\cD_{k-1}(\bar{\psi}^{2})$, respectively.
In particular, $\proj\cT_0^k(0,\bar{\psi}^{1})$ and $\proj\cT_1^k(0,\bar{\psi}^{1})$ are disjoint (cf \eqref{def:shrinking} and \eqref{def:compass}).



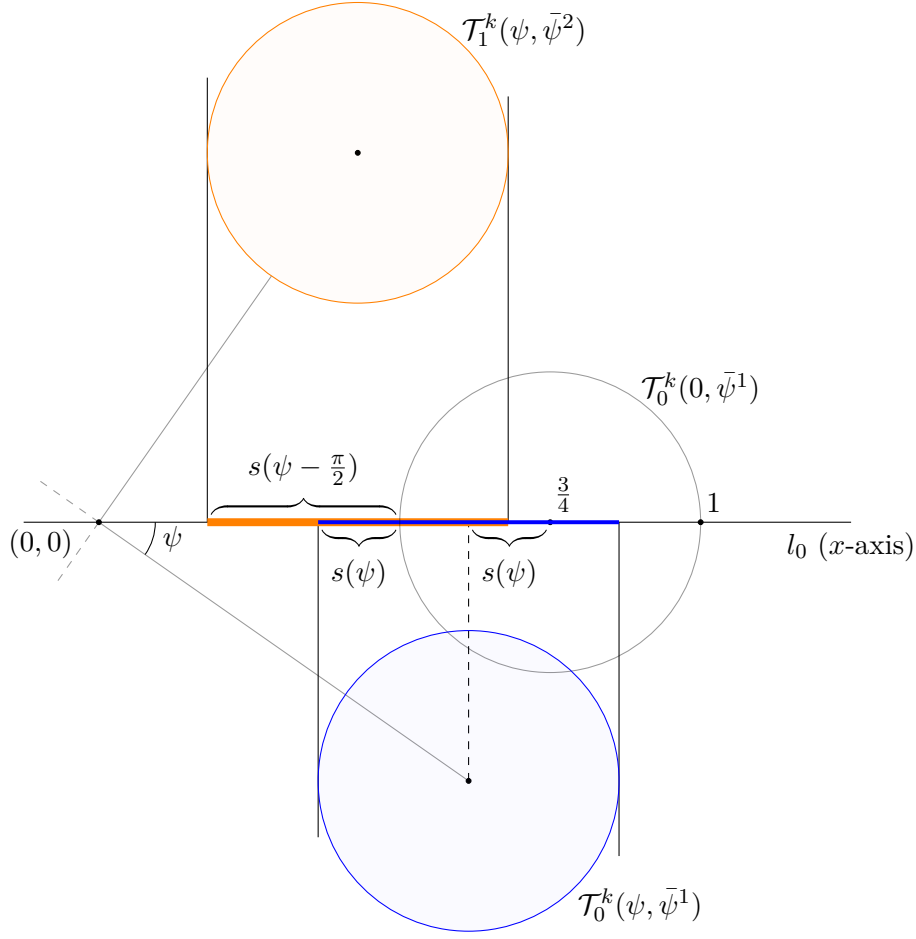
\begin{figure}
	\centering
	\begin{tikzpicture}[scale=8,decoration={pre length=1.5pt,post length=1.5pt}]
		\draw [fill] (0,0)	circle (1/250) node [below=8pt,left=6pt] at (0,0) {$(0,0)$};
		\draw [fill] ({3/4*cos(-35)},{3/4*sin(-35)}) circle (1/250) node {};
		\draw [fill] ({3/4*cos(90-35)},{3/4*sin(90-35)}) circle (1/250) node {};
		\draw [fill] (1,0)	circle (1/250) node [right=5pt,above] {$1$};
		\draw [fill] (3/4,0)circle (1/250) node [right=5pt,above] {$\frac{3}{4}$};
		\draw (-1/8,0) -- (5/4,0) node [below] {$l_0$ ($x$-axis)};
		\draw [opacity=0.4,rotate=-35] (0,0) -- (3/4,0);
		\draw [opacity=0.4,rotate=-35] (0,0) -- (0,1/2);
		\draw [opacity=0.4,rotate=-35,dashed] (0,0) -- (-1/8,0);
		\draw [opacity=0.4,rotate=-35,dashed] (0,0) -- (0,-1/8);
		\draw [dashed] ({3/4*cos(-35)},{3/4*sin(-35)}) -- ({3/4*cos(-35)},0);
		\draw ({3/4*cos(-35)-1/4},{3/4*sin(-35)-3/32}) -- ({3/4*cos(-35)-1/4},0);
		\draw ({3/4*cos(-35)+1/4},{3/4*sin(-35)-2/16}) -- ({3/4*cos(-35)+1/4},0);
		\draw ({3/4*cos(90-35)-1/4},{3/4*sin(90-35)+2/16}) -- ({3/4*cos(90-35)-1/4},0);
		\draw ({3/4*cos(90-35)+1/4},{3/4*sin(90-35)+3/32}) -- ({3/4*cos(90-35)+1/4},0);
		\draw [orange,line width=3pt] ({3/4*cos(90-35)-1/4},0) -- ({3/4*cos(90-35)+1/4},0);
		\draw [blue,line width=1.5pt] ({3/4*cos(-35)-1/4},0) -- ({3/4*cos(-35)+1/4},0);
		
		\draw [opacity=0.4] (3/4,0) circle (1/4);
		\draw [fill,blue,fill opacity=0.02] ({3/4*cos(-35)},{3/4*sin(-35)}) circle (1/4);
		\draw [fill,orange,fill opacity=0.02] ({3/4*cos(90-35)},{3/4*sin(90-35)}) circle (1/4);
		\draw (1.5/16,0) arc (0:-35:1.5/16) node [midway,right] {$\psi$};
		\draw [thick,decorate,decoration={calligraphic brace,raise=3pt,amplitude=6pt}] (1/2,0) -- ({3/4*cos(-35)-1/4},0)
			node [midway,below=10pt] {$s(\psi)$};
		\draw [thick,decorate,decoration={calligraphic brace,raise=3pt,amplitude=6pt}] (3/4,0) -- ({3/4*cos(-35)},0)
			node [midway,below=10pt] {$s(\psi)$};
		\draw [thick,decorate,decoration={calligraphic brace,raise=3pt,amplitude=6pt}] ({3/4*cos(90-35)-1/4},0) -- (1/2,0)
			node [midway,above=10pt] {$s(\psi-\frac{\pi}{2})$};
		\node [below=-4pt] at (1,1/4) {$\cT_0^k(0,\bar{\psi}^1)$};
		\node [right=8pt,above] at ({3/4*cos(-35)+1/4},{3/4*sin(-35)-1/4}) {$\cT_0^k(\psi,\bar{\psi}^1)$};
		\node [left=-6pt,below] at ({3/4*cos(90-35)+1/4},{3/4*sin(90-35)+1/4}) {$\cT_1^k(\psi,\bar{\psi}^2)$};
	\end{tikzpicture}
	\caption{The projections of \textcolor{blue}{$\cT_0^k(\psi,\bar{\psi})$} and \textcolor{orange}{$\cT_1^k(\psi,\bar{\psi})$} are contained in copies of the interval $(\frac{1}{2},1)$ shifted by \textcolor{blue}{$s(\psi)$} and \textcolor{orange}{$s(\psi-\frac{\pi}{2})$}, respectively.}
	\label{fig:shifting_projections}
\end{figure}

A simple geometric consideration (see \Cref{fig:shifting_projections}) shows that the projections of $\cT_0^k(0,\bar{\psi}^1)$ and $\cT_0^k(\psi,\bar{\psi}^1)$ differ only by a shift of $s(\psi)$, i.e.
\[\proj\cT_0^k(\psi,\bar{\psi}^1)=s(\psi)+\proj\cT_0^k(0,\bar{\psi}^1).\]
Similarly, the projections of $\cT_1^k(\psi,\bar{\psi}^2)$ and $\cT_0^k(0,\bar{\psi}^1)$ differ (on average) by a shift of $s(\psi-\frac{\pi}{2})$.

As a consequence, the events
\[\{\needle(t)\cap\cT_0^k(\psi,\bar{\psi}^{1})\not=\emptyset\} \inline{and} \{\needle(t)\cap\cT_1^k(\psi,\bar{\psi}^{2})\not=\emptyset\}\]
might not have the same probability, but their probabilities are equal to
\[\P_{\bar{\psi}^1}\left(\needle(t')\cap\cT_0^k(0,\bar{\psi}^1)\right)\]
for some appropriately shifted $t'$, since $\bar{\psi}^1$ and $\bar{\psi}^2$ are independent. We make this explicit in the following lemma.
\begin{lem}	\label{lem:shift}
	With notation as above, it holds that
	\begin{gather}
		\P_{\bar{\psi}^1}\left(\needle(t)\cap\cT_0^k(\psi,\bar{\psi}^1)\not=\emptyset\right) = \P_{\bar{\psi}^1}\left(\needle(t+s(\psi))\cap\cT_0^k(0,\bar{\psi}^1)\not=\emptyset\right)	\notag\\
	\intertext{and}
	\label{eq:shift}
		\P_{\bar{\psi}^2}\left(\needle(t)\cap\cT_1^k(\psi,\bar{\psi}^2)\not=\emptyset\right) = \P_{\bar{\psi}^1}\left(\needle(t+s(\psi-\frac{\pi}{2}))\cap\cT_0^k(0,\bar{\psi}^1)\not=\emptyset\right).
	\end{gather}
\end{lem}


With \Cref{lem:shift} at hand along with \eqref{eq:projection}, we can rewrite the overlap with our current notation in a more convenient way. First, let us denote
\begin{equation}	\label{F}
	F(t):=\P_{\bar{\psi}^1}\left(\needle(t)\cap\cT_0^k(0,\bar{\psi}^1)\right).
\end{equation}
Then, we see that 
\begin{align*}
	\cE(\psi)	&	{}=\int\P_{\bar{\psi}^1}\left(\needle(t)\cap\cT_0^k(\psi,\bar{\psi}^1)\not=\emptyset\right) \cdot \P_{\bar{\psi}^2}\left(\needle(t)\cap\cT_1^k(\psi,\bar{\psi}^2)\not=\emptyset\right)dt\\
			&	{}=\int \P_{\bar{\psi}^1}\left(\needle(t+s(\psi))\cap\cT_0^k(0,\bar{\psi}^1)\not=\emptyset\right) \cdot	\P_{\bar{\psi}^1}\Big(\needle(t+s(\psi-\frac{\pi}{2}))\cap\cT_0^k(0,\bar{\psi}^1)\not=\emptyset\Big)dt\\
			&	{}=\int F(t+s(\psi)) \cdot F(t+s(\psi-\frac{\pi}{2}))dt,
\end{align*}
where the first equality is simply the definition of $\cE$.

At this point, if we integrate over $\psi\in[0,\frac{\pi}{2}]$, we get that the
\[\text{Expectation of the overlap} =\int\cE(\psi)d\psi=\int\int F(t+s(\psi))\cdot F(t+s(\psi-\frac{\pi}{2}))d\psi dt.\]
Let's make this change of variables: $u=t+\frac34(1-\cos \psi)$ and $v=t+\frac34(1-\cos(\psi-\frac{\pi}{2}))$. The Jacobian of this change is at most $\frac{3\sqrt2}{4}$, and thus
\[\text{Expectation of the overlap} \geq\frac{4}{3\sqrt2}\int\int F(u)F(v)dudv=\frac{2\sqrt2}{3}\left(\int F(t)dt\right)^2.\]


Now, we can revert to our initial notation. And since there is no dependence on $\bar{\omega}_{n-k+1}^3$ or $\bar{\omega}_{n-k+1}^4$, we get
\begin{align*}
	\E_{\bar{\omega}_{n-k+1}^1} &	\left|\proj\cT_0^k(\bar{\omega}_{n-k+1}^1)\cap\proj\cT_1^k(\bar{\omega}_{n-k+1}^1)\right|\\
	&	=\text{Expectation of the overlap}\\
	&	\geq \frac{2\sqrt2}{3}\left(\E_{\bar{\omega}_{n-k+1}^1}\left|\proj\cT_0^k(\bar{\omega}_{n-k+1}^1)\right|\right)^2\\
	&	\xlongequal{\eqref{eq:zoom}}	\frac{2\sqrt2}{3}\cdot\frac{1}{16}\left(\E_{\bar{\omega}_{n-k+2}^1}\left|\proj\cD_{k-1}(\bar{\omega}_{n-k+2}^1)\right|\right)^2\\
	&	=\frac{1}{12\sqrt2}(D_{k-1}^1)^2.
\end{align*}

Finally, combining the fact that
\[\E_{\bar{\omega}_{n-k+1}^1}\left|\proj\cT_0^k(\bar{\omega}_{n-k+1}^1)\cap\proj\cT_1^k(\bar{\omega}_{n-k+1}^1)\right|\geq \frac{1}{12\sqrt2}(D_{k-1}^1)^2\]
with \eqref{eq:extra_overlap} and setting $c=12\sqrt2$ we get
\[D_k^1\leq D_{k-1}^1-c^{-1}(D_{k-1}^1)^2.\]
\Cref{lem:square_overlap} is proved.

\section{Comparison with \protect{\cite{Zha2019}} and \protect{\cite{PerSol2002}}} \label{sec:compa}

The random Cantor set in \cite{Zha2019} is a very close relative of the random Cantor set in this note, the difference is that Zhang's random construction of $n$ generations has $n$ independent rotations involved, whereas our construction has $1+\dots+4^{n-1}$ independent rotations. There the disks of generation $k$ are rotated by the same angle $\omega_k$, while in this note we have $4^{k-1}$ independent rotations of disks of generation $k$. Naturally, it is more difficult to work in a more chaotic model such as ours, and the techniques here use independence in a more involved way than in \cite{Zha2019}. It is just a little harder to make sense of the combinatorics involved in our model.

On the other hand, there are many \enquote{common places}: the use of overlap as the way to see the rate of decays of successive approximations of the random Cantor set, the use of \Cref{lem:square_overlap}, as well as the technical \Cref{lem:shift}.

\bigskip

Concerning \cite{PerSol2002}, there are two main differences which create difficulties. The first is the fact that at most two of the projections $\proj_\theta\cT_\alpha^k(\bar{\omega}_{n-k+1}^1)$ can intersect at each point on the line $l_\theta$. This is equivalent to line $\needle_\theta(t)$ intersecting at most two of the disks for any $t$, and is key to the square factor appearing in our calculations.

However, this is simply not true in the case of squares. In fact, in the Peres and Solomyak case the corresponding line $\needle_\theta(t)$ can simultaneously intersect $3$ squares of generation $k$ for any $k$ and any $t$. Because of this, the inequalities appearing here cannot be translated directly in the square setting.

\medskip

But even if this wasn't an obstacle, the reader should pay attention to \Cref{lem:shift}. Let's pretend that we can repeat everything before this lemma for the model of Peres and Solomyak. The role of the angle $\omega^1_{n-k+1}$ will be played by the \enquote{Favard angle} $\theta$, the shift function $s(\omega^1_{n-k+1})$ will be replaced by
\[s(\theta)=\frac{1}{2}\sin \theta,\]
and all seems to be following smoothly along the same lines. Also, the following equality
\begin{equation} \label{sPS}
	\int\E\mathbf{1}_{\left\{\needle_\theta(t)\cap\cT_0^k(\bar{\omega}_{n-k+1}^1)\not=\emptyset\right\}}d\theta=\int\E\mathbf{1}_{\left\{\needle_\theta(t+s(\theta))\cap\cT_1^k(\bar{\omega}_{n-k+1}^1)\not=\emptyset\right\}}d\theta,
\end{equation}
which would be the analogue of \eqref{eq:shift}, makes sense in principle if we understand $\omega$'s as the random variables in the Peres--Solomyak model, which assume the values $0,1,2,3$ (instead of values in the interval $[0,\frac{\pi}{2}]$ as in our's and Zhang's models).

But, there is a caveat. We reduced the function of two variables 
\[G(\psi,t):=\P_{\bar{\omega}_{n-k+2}^1}\left(\needle(t)\cap\cT_0^k(\psi,\bar{\omega}_{n-k+2}^1)\not=\emptyset\right)\]
to the composition with a function of one variable and the shift (see \eqref{F} for the definition of~$F$):
\begin{equation} \label{G}
	G(\psi,t)=G(0,t+s(\psi))=F(t+s(\psi))
\end{equation}
thanks to \eqref{eq:shift}. But looking at \eqref{sPS}, we can notice that the function 
\[\cG(\theta,t):=\E\mathbf{1}_{\left\{\needle_\theta(t)\cap\cT_0^k(\bar{\omega}_{n-k+2}^1)\not=\emptyset\right\}}\]
cannot be written as some $ \cF(t+S(\theta))$.

As a result of this misfortune, we cannot write
\[\text{Expectation of the overlap}=\int\cE(\theta)d\theta=\int\int\cF(t+S(\theta))\cdot\cF(t)d\theta dt\]
as before. Working similarly, this would in turn bring about the term $(\int\cF dt)^2$. Instead, we only have that
\[\text{Expectation of the overlap}=\int\cE(\theta)d\theta=\int\int\cG(\theta,t)\cdot\cG(\theta,t+S(\theta)))d\theta dt,\]
and it is not clear (at least to us) how to estimate this integral from below as no change of variables seems to be of help.

\section*{Acknowledgements}
We are grateful to the referee for their valuable feedback and comments.

\printbibliography	

\end{document}